\documentclass[11pt]{amsart}

\usepackage{amsmath}
\usepackage{amssymb}
\usepackage{mathrsfs}
\usepackage[all]{xy}
\usepackage{todonotes}

\newtheorem{theorem}{Theorem}[section]
\newtheorem{claim}[theorem]{Claim}

\newtheorem{lemma}[theorem]{Lemma}

\newtheorem{proposition}[theorem]{Proposition}
\newtheorem{corollary}[theorem]{Corollary}

\newcounter{maintheorem}
\newtheorem{matheorem}[maintheorem]{Main Theorem}
\newtheorem{theoremint}[maintheorem]{Theorem}
\theoremstyle{definition}
\newtheorem{definition}[theorem]{Definition}

\newtheorem{question}[theorem]{Question}

\theoremstyle{remark}
\newtheorem{remark}[theorem]{Remark}

\newtheorem{fact}[theorem]{Fact}

\def\l{{\langle}}
\def\r{{\rangle}}

\newcount\skewfactor
\def\mathunderaccent#1#2 {\let\theaccent#1\skewfactor#2
\mathpalette\putaccentunder}
\def\putaccentunder#1#2{\oalign{$#1#2$\crcr\hidewidth
\vbox to.2ex{\hbox{$#1\skew\skewfactor\theaccent{}$}\vss}\hidewidth}}

\newcommand{\one}{\mathop{1\hskip-3pt {\rm l}}}

\def\smallbox#1{\leavevmode\thinspace\hbox{\vrule\vtop{\vbox
   {\hrule\kern1pt\hbox{\vphantom{\tt/}\thinspace{\tt#1}\thinspace}}
   \kern1pt\hrule}\vrule}\thinspace}

\DeclareMathOperator{\ro}{\mathcal{RO}}

\DeclareMathOperator{\Succ}{Succ}

\DeclareMathOperator{\Add}{Cohen}

\newcommand{\cf}{{\rm cf}}
\newcommand{\id}{{\rm id}}
\newcommand{\crit}{{\rm crit}}

\author{Tom Benhamou}
\author{Gabriel Goldberg}
\title{On Cohen Generic Sets in Supercompact Prikry Extensions}
\date{\today}
\begin{document}
\begin{abstract}
    Let $\kappa$ be $\lambda$-supercompact for a regular $\lambda>\kappa$. We show that if $\lambda$ is not the successor of a singular of cofinality below $\kappa$, $\diamondsuit(S^\lambda_{\geq\kappa})$ implies that every ${<}\kappa$-strategically closed forcing of size at most $\lambda$ is a projection of the $\lambda$-supercompact Prikry forcing. We show that the situation for successors of singulars of small cofinality is independent. Finally, we provide more limitations about the class of $\kappa$-distributive projections of the $\lambda$-supercompact Prikry forcing and give examples of $\kappa$-distributive forcings of size $\lambda$ which cannot be projections of this forcing. 
\end{abstract}
\maketitle
\section{Introduction}
In this paper we investigate the possible forcings which yield intermediate generic extensions of a supercompact Prikry extension. This investigation falls into the theme of investigating a mathematical structure by classifying its substructures. In forcing theory, the classification problem of subforcings of a given forcing $\mathbb{P}$ reduces to classifying the forcings $\mathbb{Q}$ which yield intermediate generic extensions of a generic extension by $\mathbb{P}$. 

A full classification is known for many of the central forcings. For example, it is well-known that any subforcing of Cohen forcing is isomorphic to Cohen forcing\footnote{ Indeed, any atomless countable forcing is forcing equivalent to Cohen forcing.}. Maharam showed \cite{maharam} that any subforcing of the Random forcing is isomorphic to the Random forcing. Other examples include Sacks forcing \cite{Sacks1971} which yields a minimal generic extension, one that has no proper intermediate models.

In the realm of Prikry-type forcings, Gitik-Kanovei-Koepke \cite{PrikryCaseGitikKanKoe} showed that if $U$ is a normal ultrafilter over $\kappa$, then the only non-trivial subforcings of the Prikry forcing with $U$ are forcing equivalent to Prikry forcing with $U$. This was later extended by Benhamou-Gitik to the Magidor-Radin forcing \cite{TomMoti,partOne,Parttwo}.  Also, Koepke-R\"asch-Schlicht \cite{MinimalPrikry} showed that the tree Prikry forcing with distinct normal measures is minimal in the sense of Sacks. 

While the above examples demonstrate situations where the class of subforcings is narrow, there are forcings which are universal in the sense that they absorb every forcing from a given class. A typical example is the collapsing poset $Col(\kappa,\lambda)$ which produces a generic filter for every $\kappa$-strategically closed forcing of cardinality at most $\lambda$. Another example is obtained by Prikry forcing when giving up normality. Benhamou-Hayut-Gitik \cite{Benhamou_Hayut_Gitik} and Benhamou-Gitik \cite{BENHAMOU_GITIK_2024} showed that a wide variety of forcings can be absorbed into a non-normal tree Prikry forcing, and this class may depend on the initial large cardinal assumption. In fact, if $\kappa$ is strongly compact, then every $\kappa$-distributive forcing of cardinality at most $\kappa$ is a projection of the tree Prikry forcing for some tree of measures \cite[Claim 12]{Benhamou_Hayut_Gitik}.

More general and relevant to this paper is the strongly compact Prikry forcing (see Definition \ref{def: general Prikry forcing}) which we denote by $Pr_T(\mathcal{U})$:
\begin{theorem}[Folklore]
Suppose that $\kappa$ is $2^\lambda$-strongly compact. Then for every $\kappa$-distributive forcing $\mathbb{Q}$ of size at most $\lambda$ and for every fine $\kappa$-complete ultrafilter $\mathcal{U}$ on $P_\kappa(2^\lambda)$, there is a weak projection\footnote{See Definition \ref{def: projections}.} of $Pr_T(\mathcal{U})$ onto $\mathbb{Q}$. 
\end{theorem}
In particular, if $\kappa$ is supercompact then every $\kappa$-distributive forcing of cardinality at most $\lambda$ is a projection of the $2^\lambda$-supercompact Prikry forcing.

Recently, the $2^\lambda$-strongly compact Prikry forcing was brought down to $\lambda$-strongly compact Prikry:
\begin{theorem}[Benhamou-Thei-Weltsch]\label{distributive forcing can be a projection}
    Suppose that $\kappa$ is $2^\lambda$-strongly compact. Then for every $\kappa$-distributive forcing $\mathbb{Q}$ of size $\leq\lambda$ there is a fine $\kappa$-complete ultrafilter $\mathcal{U}$ on $P_\kappa(\lambda)$ such that $\mathbb{Q}$ is a weak projection of $Pr_T(\mathcal{U})$.
\end{theorem}
The above result raises the following questions which were asked in \cite{BenhamouTheiWeltsch2026}:
\begin{question}\ {}
    \begin{enumerate}
        \item Assuming that $\kappa$ is $2^\lambda$-supercompact, is every $\kappa$-distributive forcing of size $\leq\lambda$ a projection of the $\lambda$-supercompact Prikry forcing?
        \item What about $Add(\lambda,1)$ and $Col(\kappa,\lambda)$?
    \end{enumerate}
\end{question}
The main results of this paper answer these questions:
\begin{matheorem}\label{Main Theorem 1}
    Assume that $\lambda>\kappa$ is regular and is not the successor of a singular of cofinality less than $\kappa$ and assume $\diamondsuit(S^\lambda_{\geq\kappa})$. Then for every ${<}\kappa$-strategically closed forcing $\mathbb{Q}$ of size $\leq\lambda$ and for every normal ultrafilter $\mathcal{U}$ on $P_\kappa(\lambda)$, there is a weak projection of $Pr_T(\mathcal{U})$ to $\mathbb{Q}$.
\end{matheorem}
This gives a positive answer to question $(2)$. Theorem \ref{Main Theorem 1} stands in contrast to the case $\lambda=\kappa$, that is, the vanilla Prikry forcing with a normal ultrafilter on $\kappa$. Indeed, by Gitik-Kanovei-Koepke \cite{PrikryCaseGitikKanKoe}, if $U$ is a normal ultrafilter on $\kappa$, then $\Add(\kappa,1)$ (and many more) cannot be projections of Prikry forcing $Pr(U)$. In fact we prove Main Theorem \ref{Main Theorem 1} for a slightly larger class of forcings (see Theorem \ref{General Main Theorem 1}) which we call the ${<}\kappa$-diagonally distributive forcings (see Definition \ref{def filtration dist}). We then show in Proposition \ref{Proposition in btween class} that
$${<}\kappa\text{-strategically closed}\Longrightarrow {<}\kappa\text{-diagonally distributive}\Longrightarrow \kappa\text{-distributive}.$$
We also show that the above implications cannot be reversed.

The case for successors of singulars of cofinality less than $\kappa$ is more delicate. We are able to show the following independence result:
\begin{theoremint}
Relative to infinitely many supercompact cardinals, it is independent whether $Add(\kappa^{+\omega+1},1)$ is a projection of the $\kappa^{+\omega+1}$-supercompact Prikry forcing.    
\end{theoremint}

As for $(1)$, we obtain the following negative answer:
\begin{matheorem}\label{Main Theorem 2}
    Assume that $\lambda=\kappa^+$. There is a $\kappa$-distributive forcing of size $\leq\lambda$ which cannot be the projection of the $\lambda$-supercompact Prikry forcing.
\end{matheorem}
The forcing in Main Theorem \ref{Main Theorem 2} is the forcing for shooting a club through fat stationary sets. 
\section{Preliminaries}
\subsection{Forcing preliminaries} Let us start with the notion of projection:
\begin{definition}\label{def: projections}
        Let $\mathbb{P}$ and $\mathbb{Q}$ be forcing posets. 
         A \emph{weak projection} is a map $\pi\colon \mathbb{P}\rightarrow\mathbb{Q}$ such that: \begin{enumerate}
            \item For all $p,p'\in\mathbb{P}$, if $p\geq p'$ then $\pi(p)\geq \pi(p')$;
            \item for all $p\in\mathbb{P}$, there is $q^*\leq \pi(p)$ such that for every  $q\leq q^*$, there is $p'\leq p$ such that $\pi(p')\leq q$. 
        \end{enumerate}  
        A projection is a weak projection such that:
        \begin{enumerate}
            \item $\pi(\one_{\mathbb{P}})=\one_{\mathbb{Q}}$;
            \item [$(2^*)$] in condition $(2)$ above we can take $q^*=\pi(p)$. (Namely, $\pi[\mathbb{P}/p]$ is dense in $\mathbb{Q}/\pi(p)$).
        \end{enumerate}
    \end{definition}
\begin{remark}
        Condition $(2)$ is equivalent to the requirement that for every $p\in\mathbb{P}$ and for every  $D\subseteq\mathbb{Q}$ dense above $\pi(p)$,  $\pi^{-1}[D]$ is dense in $\mathbb{P}$ above $p$.
    \end{remark}

    \begin{proposition} Let $\mathbb{P},\mathbb{Q}$ be forcing notions. \begin{enumerate}
        \item The following are equivalent:
        
    \begin{enumerate}
        \item  There is a weak projection $\pi:\mathbb{P}\to RO(\mathbb{Q})$.
        \item For every $V$-generic $G\subseteq\mathbb{P}$ there is a $V$-generic $H\subseteq\mathbb{Q}$ such that $H\in V[G]$.
    \end{enumerate}
    \item The following are equivalent:
    \begin{enumerate}
    \item  There is a projection $\pi:\mathbb{P}\to RO(\mathbb{Q})$.
        \item There is a $\mathbb{P}$-name $\dot{H}$ such that 
        $\one_{\mathbb{P}}\Vdash`` \dot{H}\text{ is }V\text{-generic for }\mathbb{Q}$", and for every $V$-generic $H\subseteq \mathbb{Q}$ there is a $V$-generic $G\subseteq\mathbb{P}$ such that $(\dot{H})_G=H$.
    \end{enumerate}
    \end{enumerate}
    \end{proposition}

     The above proposition suggests that in order to analyze (weak) projections of $\mathbb{P}$, it suffices to look at intermediate models $V\subseteq M\subseteq V[G]$ of a generic extension by $\mathbb{P}$ and study forcings $\mathbb{Q}$ such that $M$ can be a generic extension by $\mathbb{Q}$. 
     
     We will consider the following combinatorial properties of a forcing notion:
\begin{definition}
    Let $\mathbb{P}$ be a forcing notion, and let $\kappa$ be regular. We say that:
    \begin{enumerate}
        \item $\mathbb{P}$ is $\kappa$-directed closed if every $X\in [\mathbb{P}]^{<\kappa}$ consisting of pairwise compatible elements is bounded in $\mathbb{P}$.
        \item $\mathbb{P}$ is $\kappa$-closed if for any $\lambda<\kappa$ and any $\langle p_\alpha\mid \alpha<\lambda\rangle\subseteq \mathbb{P}$ decreasing there is a lower bound.
        \item $\mathbb{P}$ is $\kappa$-strategically closed if Good has a winning strategy in the game $G_\kappa(\mathbb{P})$ defined as follows:
        Let $\gamma$ be an ordinal. We define the game $G_\gamma(\mathbb{P})$ between two players Bad and Good. Before round $0$, Bad chooses a starting condition $p\in\mathbb{P}$, and every subsequent move lies in $\mathbb{P}/p$. At a round $\rho<\gamma$, a sequence $\langle p_i,q_i\mid i<\rho\rangle$ has been played and Good first plays $p_\rho$, a lower bound for all previously played conditions, and then Bad plays $q_\rho\leq p_\rho$.
         Good wins iff both players can play every stage $\rho<\gamma$.
        \item $\mathbb{P}$ is ${<}\kappa$-strategically closed if for every $\rho<\kappa$, Good has a winning strategy in the game $G_\rho(\mathbb{P})$.
        \item $\mathbb{P}$ is $\kappa$-distributive if the intersection of fewer than $\kappa$ many dense open sets is dense open.
    \end{enumerate}
\end{definition}
\begin{fact} \ {}
\begin{enumerate}
    \item $\kappa\text{-dir. closed}\Rightarrow\kappa\text{-closed}\Rightarrow\kappa\text{-str.closed}\Rightarrow{<}\kappa\text{-str.closed}\Rightarrow\kappa\text{-dist.}$
    \item $\mathbb{Q}$ is $\kappa$-dist. iff $\mathbb{Q}$ does not add $<\kappa$-sequences of ordinals.
\end{enumerate}
\end{fact}
The following is due to Dobrinen:
\begin{proposition}[{\cite[Corollary 1.6]{Dobrinen2003}}]
    Let $\kappa$ be a limit cardinal. Then the following are equivalent for every forcing notion $\mathbb{Q}$:
    \begin{enumerate}
        \item $\mathbb{Q}$ is $\kappa$-distributive.
        \item Bad does not have a winning strategy in the games of length ${<}\kappa$.
    \end{enumerate}
\end{proposition}
Observe that $\mathbb{Q}$ being $\kappa$-distributive is equivalent to saying that for every $q$, the filter $\mathcal{D}_q(\mathbb{Q})$ generated by $$\{D\subseteq \mathbb{Q}\mid D\text{ is dense open above }q\}$$
is $\kappa$-complete.

\subsection{Ultrafilter and Ultrapower Preliminaries}
We denote by $P_\kappa(A)$ the collection of all subsets of $A$ of cardinality strictly less than $\kappa$. We now establish the basic terminology to be used in this paper. For more information about the ultrafilter theory on $P_\kappa(\lambda)$ we refer to \cite{Jech2003} and \cite{GoldbergUA}. \begin{definition} An ultrafilter $\mathcal{U}$  over $P_\kappa(\lambda)$ is called:
 \begin{enumerate}
     \item fine if for every $\alpha<\lambda$, $\{x\mid \alpha\in x\}\in \mathcal{U}$;
     \item normal, if it is a fine, $\kappa$-complete ultrafilter closed under diagonal intersection, namely if $(A_\alpha\mid \alpha<\lambda)\subseteq \mathcal{U}$, then 
     $$\Delta_{\alpha<\lambda}A_\alpha:=\{x\in P_\kappa(\lambda)\mid \forall\alpha\in x, \ x\in A_\alpha\}\in \mathcal{U}.$$
 \end{enumerate}
     
 \end{definition}
Given a model of set theory $M$, and an $M$-ultrafilter $U$ over $X\in M$ \footnote{i.e. $(M,U)\models U\text{ is an ultrafilter over }X$.}, we will denote the (transitive collapse of the) ultrapower of $M$ by $U$ as $M_U$ and the associated map as $j^M_U : M \to M_U$ given by $j_U(x)=[c_x]_U$. If $M=V$ is the universe, we denote $j^V_U=j_U$.
 \begin{proposition}
     Suppose that $\mathcal{U}$ is a $\kappa$-complete ultrafilter over $P_\kappa(\lambda)$. Then:
     \begin{enumerate}
         \item $\mathcal{U}$ is fine iff $j_{\mathcal{U}}''\lambda\subseteq [id]_{\mathcal{U}}$.
         \item $\mathcal{U}$ is normal iff $j_{\mathcal{U}}''\lambda=[id]_{\mathcal{U}}$.
\end{enumerate}
For the proof of the following propositions see \cite{GoldbergUA}.

 \end{proposition}
 \begin{proposition}
     Let $U$ be an $M$-ultrafilter and let $\lambda$ be a cardinal, where $M$ is closed under $\lambda$-sequences. Then the following are equivalent:
     \begin{enumerate}
         \item $(M_U)^\lambda\cap M\subseteq M_U$ (i.e. $M_U$ is closed under $\lambda$-sequences from $M$).
         \item $j_U''\lambda\in M_U$.
     \end{enumerate}
 \end{proposition}
 \begin{remark}
     We note that if $\kappa=\lambda$, then a $P_\kappa(\lambda)$ fine $\kappa$-complete measure concentrates on ordinals and therefore we may identify it with a measure on $\kappa$. Also, normality takes its usual meaning. In this sense, the ultrafilter theory for $P_\kappa(\lambda)$ generalizes that for ultrafilters over $\kappa$.
 \end{remark}
 In the context of large cardinals, the most common way of producing ultrafilters is to let $j : M \to N$ be an elementary embedding, $X\in M$ any set, and $x \in j(X) \cap N$.
    We derive an $M$-ultrafilter $D(x,j)$ over $X$ as follows: for any $A \subseteq X$ in $M$,
    \[ A \in D(x,j) \iff x \in j(A) \]
By considering $j_{D(x,j)}:M\to M_{D(x,j)}$, we can complete the diagram using the map $k : M_{D(x,j)} \to N$ defined by $k([f]_{D(x,j)}) = j(f)(x)$. This is an elementary embedding such that $k \circ j_{D(x,j)} = j$
and $k([\text{id}]_{D(x,j)}) = x$.
    
If we started with $j=j_U$ for some ultrafilter $U$ over $Y$, then the function $f:Y\to X$ with $[f]_U = x$ can be used to define the ultrafilter $D(x,j)$ as
$$f_*(U):=\{A\subseteq X\mid f^{-1}[A]\in U\}.$$
This gives rise to the Rudin-Keisler order:
 \begin{definition}[Rudin-Keisler and Kat\v{e}tov Order]\label{def rk}
    Let $U,W$ be ultrafilters over sets $X,Y$ respectively.
    We define the \textit{Rudin-Keisler order} by $U \leq_{RK} W$ iff there is a function $f : Y \to X$ such that $A \in U$ iff $f^{-1}[A] \in W$.
    We write $U = f_* W$. 
    
    Given a filter $F$ over $X$,
    we define the \textit{Kat\v{e}tov order} by $F \leq_K W$, if there is an ultrafilter $F\subseteq U$ on $X$ such that $U\leq_{RK} W$.
\end{definition}
Summarizing the above we get
\begin{proposition}
    Let $U,W$ be $M$-ultrafilters on $X,Y$ respectively. The following are equivalent:
    \begin{enumerate}
        \item $U \leq_{RK} W$.
        \item There is $f:Y\to X$ such that $U=D([f]_W,j_W)$.
        \item There is an elementary embedding $k:M_U\to M_W$ such that $k\circ j_U=j_W$.
    \end{enumerate}
\end{proposition}

\begin{corollary}
    Let $F$ be a filter over $X$ and $W$ any ultrafilter.
    Then $F \leq_K W$ iff there is some  $x \in  j_W(X)$ such that $F \subseteq D(x,j_W)$.
\end{corollary}
Let us introduce our notations for the iterated ultrapower of $V$ by an ultrafilter $U$ of length $\delta$. This is a system $\{ M^U_\alpha, j^U_{\beta,\gamma} : \alpha < \delta \text{ and } \beta < \gamma < \delta \}$
where:
\begin{enumerate}
    \item $M^U_0 = V$.
    \item $M^U_{\alpha+1}$ is the ultrapower of $M^U_\alpha$ by $j^U_{0,\alpha}(U)$ and $j^U_{\alpha,\alpha+1}$ is the associated ultrapower map.
    \item $j^U_{\xi,\alpha+1} = j^U_{\alpha,\alpha+1} \circ j^U_{\xi,\alpha}$ for $\xi < \alpha$.
    \item Whenever $\lambda$ is limit, $M^U_\lambda$ is the direct limit of $\langle M^U_\alpha : \alpha < \lambda \rangle$ via the $j^U_{\alpha,\beta}$'s for $\alpha < \beta < \lambda$, and $j^U_{\alpha,\lambda}$'s are the direct limit maps.
\end{enumerate}
We remove the superscript $U$ when there is no ambiguity.
\subsection{The strongly compact and supercompact Prikry forcing}
Let us introduce the main forcing of interest to this paper. The forcing is due to Magidor \cite{Magidor1977} and the proofs can be found in \cite{Gitik2010}. Fix a cardinal $\kappa$ which is $\lambda$-strongly compact for some cardinal $\lambda\geq\kappa$.

Let $\mathcal{U}$ be a fine $\kappa$-complete ultrafilter on $P_\kappa(\lambda)$. The Prikry forcing with $\mathcal{U}$ will produce an $\omega$-sequence in $P_\kappa(\lambda)$ with respect to the Magidor order defined by: $$x\prec y\text{ iff }x\subseteq y\text{ and }|x|<|y\cap \kappa|.$$
Denote by $[P_\kappa(\lambda)]^{<\omega}_{\prec}$ the collection of all finite $\prec$-increasing sequences from $P_\kappa(\lambda)$.
    \begin{definition}[Tree-Prikry forcing]\label{def: general Prikry forcing}
     Let $Pr_T(\mathcal{U})$ consist of conditions of the form $\langle x_0,...,x_{n-1}, T\rangle$ such that:  \begin{enumerate}
    \item $\langle x_0,x_1, ...,x_{n-1}\r\in [P_\kappa(\lambda)]^{<\omega}_\prec$;
    \item $T\subseteq [P_\kappa(\lambda)]^{<\omega}_\prec$ and $(T,\sqsubseteq)$ is a tree\footnote{Namely, it is ordered by the end-extension order $\sqsubseteq$ and closed under restrictions.};
        \item $\langle x_0,x_1, ...,x_{n-1}\r$ is the stem of $T$, that is, every $t\in T$ is $\sqsubseteq$-comparable with $\langle x_0,x_1, ...,x_{n-1}\r$;
        \item $T$ is $\mathcal{U}$-fat, namely, for each $t\in T/\langle x_0,...,x_{n-1}\rangle$ $$\text{Succ}_T(t)=\{x\in P_\kappa(\lambda)\mid t^\smallfrown x\in T\}\in \mathcal{U}.$$
    \end{enumerate} 
The order is defined by $\langle x_0,...,x_{n-1},T\rangle\leq \langle y_0,...,y_{m-1},S\rangle$ if and only if
\begin{enumerate}
    \item $\langle y_0,...,y_{m-1}\rangle\sqsubseteq \langle x_0,...,x_{n-1}\rangle$.
    \item $\langle x_0,...,x_{n-1}\rangle\in S$.
    \item $T\subseteq S$.
\end{enumerate}
We also define the direct extension order by $\langle x_0,...,x_{n-1},T\rangle\leq^* \langle y_0,...,y_{m-1},S\rangle$ if and only if $\langle x_0,...,x_{n-1},T\rangle\leq \langle y_0,...,y_{m-1},S\rangle$ and $n=m$. \qed
\end{definition}
Note that $\langle x_0,...,x_{n-1},T\rangle\leq \langle y_0,...,y_{m-1},S\rangle$ iff $T\subseteq S$. 

If $\mathcal{U}$ is normal, then we can replace the tree $T$ by a single set $A\in \mathcal{U}$:
\begin{claim}
    If $T$ is $\mathcal{U}$-fat with a stem $t_0$, and $\mathcal{U}$ is normal, then there is $A\in\mathcal{U}$ such that $$t_0^\smallfrown[A]^{<\omega}_\prec:=\{t_0^\smallfrown s\mid s\in [A]^{<\omega}_\prec\} \subseteq T.$$
\end{claim}
\begin{proof}
    Set
    $A=\{x\mid t_0\prec x\}\cap \Delta_{t\in T/t_0}\Succ_T(t)$ where $$\Delta_{t\in T/t_0}\Succ_T(t):=\{x\in P_\kappa(\lambda)\mid \forall s\in T,s\prec x\Rightarrow x\in \Succ_T(s)\}.$$ Then $A\in \mathcal{U}$. Indeed, if $s\in j_{\mathcal{U}}(T)$ and $s\prec j_{\mathcal{U}}''\lambda$, then there is $s'\in T$, such that $j_{\mathcal{U}}(s')=s$. Since $\Succ_T(s')\in \mathcal{U}$, it follows that $$j_{\mathcal{U}}''\lambda\in j_{\mathcal{U}}(\Succ_T(s'))=\Succ_{j_{\mathcal{U}}(T)}(s).$$
    Hence $j_{\mathcal{U}}''\lambda\in j_{\mathcal{U}}(A)$. To see that $t_0^\smallfrown[A]^{<\omega}_\prec\subseteq T$, take any $(x_0,...,x_{n-1})\in [A]^{<\omega}_\prec$, then $t_0\prec x_0$, and therefore $x_0\in \Succ_T(t_0)$, namely $t_0^\smallfrown x_0\in T$. Then $t_0^\smallfrown x_0\prec x_1$ and therefore $x_1\in \Succ_T(t_0\smallfrown x_0)$ and so on.
\end{proof}
We conclude the following:
\begin{corollary}
    Suppose that $\mathcal{U}$ is a normal measure on $P_\kappa(\lambda)$. Then $Pr_T(\mathcal{U})$ is forcing equivalent to the forcing $Pr(\mathcal{U})$ consisting of $(\vec{x},A)\in [P_\kappa(\lambda)]^{<\omega}_\prec\times \mathcal{U}$ and $(\vec{x},A)\leq (\vec{y},B)$ iff $\vec{y}\sqsubseteq \vec{x}$, $\vec{x}\setminus \vec{y}\subseteq B$ and $A\subseteq B$.
\end{corollary}
For a condition $p=(t,T)$, let us denote by $t(p)=t$ and $T(p)=T$. Also, for every $n<\omega$ we define $\mathcal{L}_n(T)=T\cap [P_\kappa(\lambda)]^n$.
\begin{proposition}[Properties of $Pr_T(\mathcal{U})$]
    \ {}
    \begin{enumerate}
        \item $Pr_T(\mathcal{U})$ is $(\lambda^{<\kappa})^{+}$-cc.
        \item $(Pr_T(\mathcal{U}),\leq^*)$ is $\kappa$-directed closed.
    \item The Prikry property: for every sentence $\sigma$ in the forcing language and every $p\in Pr_T(\mathcal{U})$ there is $p^*\leq^* p$ such that $p^* ||\sigma$.
    \item The strong Prikry property: for every dense open set $D$ and every $p\in Pr_T(\mathcal{U})$, there is $n<\omega$ and $p^*\leq^* p$ such that for every $s\in \mathcal{L}_n(T(p^*))$, $p^{*\smallfrown}s\in D$.
    \end{enumerate}
\end{proposition}
As a corollary we see that $Pr_T(\mathcal{U})$ does not add bounded subsets to $\kappa$, and preserves cardinals in the range $[0,\kappa]\cup [(\lambda^{<\kappa})^+,\infty)$. 
\begin{theorem}
  $Pr_T(\mathcal{U})$ collapses every cardinal in the interval $(\kappa,\lambda^{<\kappa}]$ and singularizes $\kappa$ to have cofinality $\omega$.   
\end{theorem}
We will make frequent use of the Bukovsky-Dehornoy analysis of Prikry forcing and iterated ultrapowers:
\begin{theorem}\label{Bukovsky-Dehornoy}
    Let $\mathcal{U}$ be a fine $\kappa$-complete ultrafilter over $P_\kappa(\lambda)$.
    Let $\{ M_n , j_{m,n} : m < n \in \omega\}$ be the iterated ultrapower of $V$ by $\mathcal{U}$ of length $\omega$.
    Consider the sequence
    \[ S = \{ j_{n,\omega}([\id]_{\mathcal{U}^n}) : 1\leq n < \omega \} \]
    Then $S$ is a $j_{0,\omega}(Pr_T(\mathcal{U}))$-generic sequence over $M_\omega$.
    Furthermore, $$\bigcap_{n<\omega}M_n=M_{\omega}[S].$$
    It follows that $M_{\omega}[S]$ is closed under $\kappa$-sequences from $V$.
\end{theorem}

Coming back to the question of intermediate models of the strongly compact Prikry forcing, the following \cite{BenhamouTheiWeltsch2026} provides a useful characterization inside the class of the $\kappa$-distributive forcings of size $\leq\lambda$: 
\begin{theorem}[Benhamou-Thei-Weltsch]\label{Projection criteria}
Let $\kappa$ be a $\lambda$-strongly compact cardinal. Let $\mathbb{Q}$ be a forcing of size $\leq\lambda$, and let $\mathcal{U}$ be a fine $\kappa$-complete ultrafilter over $P_\kappa(\lambda)$. Then the following are equivalent:
\begin{enumerate}
    \item $\mathbb{Q}$ is a $\kappa$-distributive weak projection with a dense image of $Pr_T(\mathcal{U})$.
    \item For every $q\in \mathbb{Q}$, there is $n_q$ such that $$\bigcap j_{\mathcal{U}^{n_q}}{}``\mathcal{D}_q(\mathbb{Q})\neq\emptyset.$$
    \item For every $q\in\mathbb{Q}$, there is $n_q$ such that $\mathcal{D}_q(\mathbb{Q})\leq_K \mathcal{U}^{n_q}$.
\end{enumerate}
\end{theorem}
\section{Constructing generic filters over ultrapowers}
\subsection{Woodin's Construction} The starting point for Main Theorem \ref{Main Theorem 1} is the following construction of Woodin: 
\begin{proposition}[Woodin]\label{Proposition: woodin's construction}
Assume $\diamondsuit(S^{\kappa^+}_\kappa)$, and let $U$ be a nonprincipal $\kappa$-complete ultrafilter over $\kappa$. Let $\mathbb{Q}\in M_U$ be a forcing which, in $M_U$, is $\kappa^+$-closed and has cardinality at most $j_U(\kappa^+)$. Then for any $q\in\mathbb{Q}$ there is an $M_U$-generic filter $G\in V$ for $\mathbb{Q}$ with $q\in G$.
\end{proposition}
\begin{proof}
Assume that $\mathbb{Q}=[\nu\mapsto\mathbb{Q}_\nu]_U$ is the ultrapower representation and suppose without loss of generality that for each $\nu$, $\mathbb{Q}_\nu\subseteq\kappa^+$. Let $\l A_\alpha\mid \alpha\in S^{\kappa^+}_\kappa\r$ be a diamond sequence. Fix a bijection $\varphi:\kappa^+\times \kappa\leftrightarrow \kappa^+$ and let $C_0$ be the club of points $\alpha$ such that $\varphi\restriction \alpha\times\kappa:\alpha\times\kappa\leftrightarrow \alpha$. For each $\alpha\in S_\kappa^{\kappa^+}$, we let $A^*_\alpha=\varphi^{-1}[A_\alpha]$. For $X\subseteq\kappa^+\times\kappa$ and $\nu<\kappa$ denote by $X_\nu=\{\alpha<\kappa^+\mid (\alpha,\nu)\in X\}$. Then the sequence $\l A^*_\alpha\mid \alpha\in S^{\kappa^+}_\kappa\r$ has the following property:

    For every sequence $(X_\nu\mid\nu<\kappa)\subseteq P(\kappa^+)$, \begin{equation}\{\alpha\in S^{\kappa^+}_\kappa\mid \langle(A^*_\alpha)_\nu\mid \nu<\kappa\rangle=\langle X_\nu\cap\alpha\mid \nu<\kappa\rangle\} \text{ is stationary in }\kappa^+.
\end{equation}
Let $\mathbb{Q}_{\alpha,\nu}=\mathbb{Q}_\nu\cap\alpha$. Let us construct a descending sequence $\l q_\xi\mid \xi<\kappa^+\r\subseteq\mathbb{Q}$ below $q$ such that for every dense set $D\subseteq\mathbb{Q}$ in $M_U$, there is $\xi<\kappa^+$ such that $q_\xi\in D$.
We will also fix (arbitrarily) at stage \(\xi\) of the construction a function $f_\xi:\kappa\to\kappa^+$ such that $f_\xi(\nu)\in\mathbb{Q}_\nu$ for each $\nu<\kappa$ and $[f_\xi]_U=q_\xi$.

Set $q_0=q$. At successor steps we set $q_{\alpha+1}=q_\alpha$, and at stages of cofinality less than $\kappa$, we set $q_\alpha$ to be a lower bound for  $\{ q_\beta\mid \beta < \alpha\}$, which exists by the $\kappa^+$-closure of $\mathbb{Q}$ in $M_U$ and the closure of $M_U$ under $\kappa$-sequences.

Suppose that $cf(\alpha)=\kappa$. We say 
that \textit{$\alpha$ is an active stage} if:
\begin{enumerate}
    \item Assumption 1: for every $\beta<\alpha$ and $\nu<\kappa$, $f_\beta(\nu)\in\mathbb{Q}_{\alpha,\nu}$.  
    \item Assumption 2: $A^*_\alpha$ codes a sequence $(E_\nu\mid \nu<\kappa)$ with $E_\nu$ dense in $\mathbb{Q}_{\alpha,\nu}$.
    \item Assumption 3: $\mathbb{Q}\cap j_U(\alpha)$ is $\kappa^+$-closed in $M_U$.
\end{enumerate}
Assumption 1 ensures that for each $\beta<\alpha$, $q_\beta = j_U(f_\beta)([\text{id}]_{U})\in \mathbb{Q}\cap j_U(\alpha)$. By Assumption 2,
 $E^*:=[\nu\mapsto E_\nu]_U$ is dense in $\mathbb{Q}\cap j_U(\alpha)$, and by Assumption 3, there is some $q_\alpha\in E^*$ such that $q_\alpha\leq q_\beta$ for every $\beta<\alpha$.

 If $cf(\alpha)=\kappa$ and $\alpha$ is not an active stage, take $q_\alpha$ to be a lower bound for $\{q_\beta\mid\beta<\alpha\}$, which exists by the $\kappa^+$-closure of $\mathbb{Q}$ in $M_U$.

This concludes the definition of $\l q_\alpha\mid\alpha<\kappa^+\r$. Now let $D\in M_U$ be a dense open subset of $\mathbb{Q}$. Write $D=[\nu\mapsto D_\nu]_U$, where each $D_\nu$ is dense open in $\mathbb{Q}_\nu$.
For each $\nu<\kappa$, there is a club $C_\nu\subseteq\kappa^+$ such that for every $\alpha\in C_\nu$, $D_\nu\cap\alpha$ is dense in $\mathbb{Q}_{\alpha,\nu}$. (Take $C_\nu$ to be the club of closure points of a function sending each $\zeta\in\mathbb{Q}_\nu$ to some $\zeta^\prime\in D_\nu$ with $\zeta^\prime\leq_{\mathbb{Q}_\nu}\zeta$.)
Let $C^*=\bigcap_{\nu<\kappa}C_\nu$. Then for every $\alpha\in C^*$ and every $\nu<\kappa$, $D_\nu\cap\alpha$ is dense in $\mathbb{Q}_{\alpha,\nu}$.
It is easy to see that Assumption $1$ also occurs on a club. Choose $\theta_\nu<\kappa$ such that $[\nu\mapsto\theta_\nu]_U=\kappa^+$. By {\L}o{\'s}'s theorem, for $U$-almost every $\nu$, $\mathbb{Q}_\nu$ is $\theta_\nu$-closed. For each such $\nu$, fix a choice of lower bounds for descending sequences of length $<\theta_\nu$. Since $\kappa^{<\kappa}=\kappa$, there is a club of $\alpha<\kappa^+$ such that whenever $\operatorname{cf}(\alpha)=\kappa$, every such sequence in $\mathbb{Q}_{\alpha,\nu}$ has its chosen lower bound below $\alpha$, simultaneously for all such $\nu$. Hence $\mathbb{Q}\cap j_U(\alpha)$ is $\kappa^+$-closed in $M_U$ at these stages. 
By the properties of \(\langle A^*_\alpha \mid \alpha \in S^{\kappa^+}_\kappa\rangle\), there is some $\alpha\in C^*\cap S^{\kappa^+}_\kappa$ satisfying Assumptions 1 and 3 such that $(D_\nu\cap\alpha)_{\nu<\kappa}$ is coded by $A^*_\alpha$. Therefore $\alpha$ is an active stage. By our construction,
\[
    q_\alpha\in[\nu\mapsto D_\nu\cap\alpha]_U=D\cap j_U(\alpha),
\]
as desired.
\end{proof}
We can now deduce a part of Main Theorem \ref{Main Theorem 1}.
\begin{theorem}\label{theorem: there is a projeciton}
    Suppose that $\diamondsuit(S^{\kappa^+}_\kappa)$ holds and let $\mathcal{U}$ be a \(\kappa\)-complete normal fine ultrafilter on $P_\kappa(\kappa^+)$. Then for every $\kappa$-closed forcing $\mathbb{P}$ of cardinality $\kappa^+$, and every $p\in\mathbb{P}$,  $\bigcap_{D\in \mathcal{D}_p(\mathbb{P})} j_{\mathcal{U}}(D)\neq\emptyset$.
\end{theorem}
\begin{proof}
    Let $U$ be the projection of \(\mathcal U\) to a normal $\kappa$-complete ultrafilter over $\kappa$. Let $i:V\to M_U$, $j:V\to M_{\mathcal{U}}$ and $k:M_U\to M_{\mathcal{U}}$ be the usual commutative diagram. Let $p\in\mathbb{P}$. By Proposition \ref{Proposition: woodin's construction}, let $G\in V$ be $M_U$-generic for $i(\mathbb{P})$, with $i(p)\in G$ and note that $|G|=|i(\mathbb{P})|=\kappa^+$ by the diamond assumption. By the $\kappa$-closure of the forcing (or simply by the way we constructed the generic), we can find a descending sequence \(\langle q_\alpha : \alpha < \kappa^+\rangle\subseteq G\) that generates \(G\). By the normality of $\mathcal{U}$, $M_{\mathcal{U}}$ is closed under $\kappa^+$-sequences of its elements, it follows that $\langle k(q_\alpha) : \alpha < \kappa^+\rangle\in M_{\mathcal{U}}$. Since $j(\mathbb{P})$ is $j(\kappa)$-closed, and $j(\kappa)>\kappa^+$, there is a lower bound $q\in j(\mathbb{P})$ for $\langle k(q_\alpha) : \alpha < \kappa^+\rangle$. Let us show that $q$ is as desired. Let $D\in \mathcal{D}_p(\mathbb{P})$. Then $i(D)\in M_U$ is dense in $i(\mathbb{P})$ above $i(p)$, and by $M_U$-genericity of $G$ and since $i(p)\in G$, there is $p'\in i(D)\cap G$. Hence $k(p')\in j(D)$, and since $D$ is open and \(q\leq k(p')\), we conclude that $q\in j(D)$, as desired.
\end{proof}
By Theorem \ref{theorem: there is a projeciton} and Theorem \ref{Projection criteria} we get:
\begin{corollary}
    Assume $\diamondsuit(S^{\kappa^+}_\kappa)$, then for every $\kappa$-closed forcing $\mathbb{P}$ of cardinality $\kappa^+$ and every normal measure $\mathcal{U}$ on $P_\kappa(\kappa^+)$, there is a weak projection from $Pr(\mathcal{U})$ onto $\mathbb{P}$.
\end{corollary}
Normality seemed crucial here, so the following question is unclear:
\begin{question}
    Is every $\kappa$-closed forcing of size $\leq\kappa^+$ (or just $Add(\kappa^+,1)$) a projection of $Pr_T(\mathcal{U})$, for every fine $\kappa$-complete ultrafilter $\mathcal{U}$ on $P_\kappa(\kappa^+)$?
\end{question}
\begin{remark}
In the case that $2^\kappa>\kappa^+$,  $Add(\kappa^+,1)$ collapses $2^\kappa$ to $\kappa^+$. On the other hand, without any cardinal arithmetic assumptions, the $\kappa^+$-supercompact Prikry forcing is $\kappa^{++}$-cc (since $(\kappa^+)^{<\kappa}=\kappa^+$) and therefore does not collapse $\kappa^{++}$. Hence if $2^\kappa>\kappa^+$, $Add(\kappa^+,1)$ cannot be a projection of $\kappa^+$-supercompact Prikry forcing.
\end{remark}
We do not know if the $\diamondsuit(S^{\kappa^+}_\kappa)$ assumption is necessary:
\begin{question}
    Is there a model of $GCH$ in which $Add(\kappa^+,1)$ is not a projection of the $\kappa^+$-supercompact Prikry forcing?
\end{question}
Note that $\diamondsuit(S^{\kappa^+}_\kappa)$ would necessarily fail in such a model.
\subsection{Main Theorem \ref{Main Theorem 1}}

Let us now state and prove our first general result from which Main Theorem \ref{Main Theorem 1} follows. Let us first introduce the class of forcings for which the theorem applies.
Suppose $\mathbb{P}$ is a poset, $p$ is a condition in $\mathbb{P}$, and $\kappa$ is a cardinal. Given a $\subseteq$-increasing sequence $\vec{N}=\l N_i\mid i<\gamma\r$ of elementary submodels of $H_\theta$ for some large enough $\theta$ such that:
\begin{enumerate}
    \item [(i)] $\gamma < \kappa$;
    \item [(ii)] for every $i<\gamma$, $(\bigcup_{j<i}N_j)^{<\kappa}\subseteq N_{i}$; 
    \item [(iii)] $\kappa+1\cup\{p,\mathbb{P}\}\subseteq N_0$.
\end{enumerate} 
We define $G(\vec{N},\mathbb{P}/p)$ to be the game between two players playing conditions $p_i,q_i\in \mathbb{P}\cap N_{i}$ above $p$ respectively. The game is of length $\gamma$ such that at stages $i<\gamma$, Good plays first, and chooses $p_i\in N_{i}$ stronger than every element of the sequence $\l p_j,q_j\mid j<i\r$ and then Bad plays $q_i\leq p_i$ also in $N_{i}$. If no such $p_i$ exists, the run terminates (and Bad wins the run).
Good wins the run iff they survive all stages $i<\gamma$ and there is a lower bound to the sequence $\l p_i\mid i<\gamma\r$.

A sequence satisfying $(i)$--$(iii)$ is called a $(\kappa,p,\mathbb{P})$-\textit{playful} sequence of elementary submodels of $H_\theta$.
\begin{definition}\label{def filtration dist}
A forcing notion $\mathbb{P}$ is said to be ${<}\kappa$-diagonally distributive if for all sufficiently large \(\theta\), and for all $(\kappa,p,\mathbb{P})$-playful sequences of elementary submodels of $H_\theta$, Bad does not have a winning strategy in \(G(\vec N,\mathbb P/p)\). 
    \end{definition}
    \begin{remark}
    \begin{enumerate}
        
        \item 
Shelah has similar definitions \cite[Definition A.1.14]{Shelah2003NotCollapsing}.
\item It is well known that Bad not having a winning strategy in $G(\mathbb{P}/p)$ is equivalent to $\kappa$-distributivity. Hence on the class of forcings of size $\leq\kappa$, $\kappa$-distributive=${<}\kappa$-diagonally distributive.
\item If $|\mathbb{P}|\leq\kappa$, then the model $N_0$ already includes $\mathbb{P}$ and therefore the game $G(\vec{N},\mathbb{P}/p)$ is the same as $G(\mathbb{P}/p)$. 
\item From the previous item it follows that if $|\mathbb{P}|\leq\kappa$ then $\mathbb{P}$ is ${<}\kappa$-diagonally distributive iff it is $\kappa$-distributive.

    \end{enumerate}
    \end{remark}
Main Theorem \ref{Main Theorem 1} follows from the following theorem and Proposition \ref{Proposition in btween class}:
\begin{theorem}\label{General Main Theorem 1}
    Let $\kappa < \lambda$ such that:
    \begin{enumerate}
        \item [(a)] $\lambda$ is not the successor of a singular of cofinality $<\kappa$.\footnote{Recall that by Solovay \cite{Solovay1974} (see also \cite{GoldbergUA}), if $\kappa$ is $\delta$ strongly compact and $cf(\delta)\geq\kappa$, then $\delta^{<\kappa}=\delta$.}
        \item [(b)] $\diamondsuit(S^{\lambda}_{\geq\kappa})$.
    \end{enumerate}  
    Then for every ${<}\kappa$-diagonally distributive forcing $\mathbb{P}$ of cardinality $\lambda$, and for every normal ultrafilter $\mathcal{U}$ on $P_\kappa(\lambda)$ and every $p\in\mathbb{P}$, $ \bigcap j_{\mathcal{U}}``\mathcal{D}_p(\mathbb{P})\neq\emptyset$. 
\end{theorem}
    \begin{proof}
Let $\l A_\alpha\mid \alpha\in  S^{\lambda}_{\geq\kappa}\r$ be a diamond sequence, and let us assume that $\mathbb{P}=\lambda$. Since $\mathcal{U}$ is normal, $\l A_\alpha\mid \alpha\in S^{\lambda}_{\geq\kappa}\r\in M_{\mathcal{U}}$. In $V$, let $\theta$ be sufficiently large and regular. Since $\lambda^{<\kappa}=\lambda$, choose $X\prec H_\theta$ with $|X|=\lambda$, $\lambda\cup\{\mathbb P,p\}\subseteq X$, and $X^{<\kappa}\subseteq X$. Fix a sequence $\vec{N}=\l N_i\mid i<\lambda\r$ of elementary submodels of $H_\theta$, each contained in $X$, such that:
\begin{enumerate}
    \item $\vec{N}$ is $\subseteq$-increasing and continuous at points of cofinality $\geq\kappa$.
    \item for each $i<\lambda$, $|N_i|<\lambda$.
    \item $(\bigcup_{j<i}N_j)^{<\kappa}\cup \sup(\bigcup_{j<i}N_j\cap\lambda)\subseteq N_{i}$.
    \item $p,\mathbb{P}\in N_0$ and $\kappa+1\subseteq N_0$.
    \item $\bigcup_{i<\lambda}N_i=X$.
\end{enumerate}
To see that such a sequence exists, we first note that since $\lambda$ is not the successor of a singular of cofinality $<\kappa$, if $\delta<\lambda$, then $\delta^{<\kappa}<\lambda$. Enumerate $X=\{x_i\mid i<\lambda\}$. Find $N_0\prec X$ with $x_0\in N_0$ and such that $(4)$ holds. At steps which are not limit points of cofinality $\geq\kappa$ (i.e. successor stages or limits of small cofinality), the regularity of $\lambda$, the induction hypothesis, and the previous observation, ensure that $\sup(\bigcup_{j<i}N_j\cap\lambda)<\lambda$ and $|(\bigcup_{j<i}N_j)^{<\kappa}|<\lambda$. Hence by L\"owenheim--Skolem, we find $N_{i}\prec X$ which includes $(\bigcup_{j<i}N_j)^{<\kappa}\cup \sup(\bigcup_{j<i}N_j\cap\lambda)$, with $|N_{i}|<\lambda$. At successor steps $i+1$ we also require that $x_i\in N_{i+1}$. At limit points of cofinality $\geq\kappa$, we must take the union. Thus $(1)$ is satisfied. $(2)$ remains true by regularity. $(3)$ follows by the induction hypothesis and the fact that the cofinality is $\geq\kappa$. Finally $(4),(5)$ follow from the construction.

Let $\mathbb{P}_\alpha=\mathbb{P}\cap N_\alpha=\lambda\cap N_\alpha$. 
Work in $M_{\mathcal{U}}$, then by elementarity $j_{\mathcal{U}}(\vec{N})$ has the same properties, and we may consider the  subsequence $j_{\mathcal{U}}''\vec{N}=j_{\mathcal{U}}(\vec{N})\restriction j_{\mathcal{U}}''\lambda$. Note that this is a $(j_{\mathcal{U}}(\kappa),{j_{\mathcal{U}}(p)},j_{\mathcal{U}}(\mathbb{P}))$-playful, $\subseteq$-increasing (nowhere continuous) sequence in $M_{\mathcal{U}}$: $(i)$ follows since this is a sequence of length $\lambda<j_{\mathcal{U}}(\kappa)$, $(ii)$ holds since for every $i<\lambda$, by elementarity we have
\[
    \left(\bigcup_{j<i}j_{\mathcal{U}}(N_j)\right)^{<j_{\mathcal{U}}(\kappa)}
    \subseteq \left(j_{\mathcal{U}}\!\left(\bigcup_{j<i}N_j\right)\right)^{<j_{\mathcal{U}}(\kappa)}
    \subseteq j_{\mathcal{U}}(N_i).
\]
$(iii)$ follows from $(4)$ and elementarity. Since $M_\mathcal{U}\models j_{\mathcal{U}}(\mathbb{P})$ is ${<}j_\mathcal{U}(\kappa)$-diagonally distributive, Bad has no winning strategy in $M_{\mathcal{U}}$ for the game $G(j_{\mathcal{U}}''\vec{N},j_{\mathcal{U}}(\mathbb{P}/p))^{M_{\mathcal{U}}}$.

Let us define a strategy for Bad in $G(j_{\mathcal{U}}``\vec{N},j_{\mathcal{U}}(\mathbb{P}/p))^{M_{\mathcal{U}}}$
as follows. Suppose that at stage $\alpha<\lambda$, we have recursively obtained a run $\l p_i,q_i\mid i<\alpha\r$ which is played according to the strategy, and Good plays $p_\alpha\in j_{\mathcal{U}}(N_{\alpha})$. 

We say that $\alpha$ is an \textit{active} stage if \(\alpha\) is a limit ordinal of cofinality at least \(\kappa\) and the following holds:

\begin{enumerate}
    \item [$(*)$]  the diamond sequence $A_\alpha$ is a dense open subset $E\subseteq \mathbb{P}_\alpha$ above $p$.
\end{enumerate}
By elementarity, $E^*:=j_{\mathcal{U}}(E)$ is dense in $j_{\mathcal{U}}(\mathbb{P}_\alpha)$, and so Bad can play some $q_\alpha\leq p_\alpha$ in $E^*$. Note that this is a legitimate move inside the model $j_{\mathcal{U}}(N_\alpha)$.

If \(\alpha\) is not an active stage, and Good played \(p_\alpha\in j_{\mathcal{U}}(\mathbb P_\alpha)\), then Bad plays \(q_\alpha = p_\alpha\).
 
This concludes the definition of the strategy for Bad. Let $\l p_\alpha,q_\alpha\mid \alpha<\lambda\r$ be a winning run in \(G(j_\mathcal U``\vec N, j_\mathcal U(\mathbb P/p))^{M_{\mathcal{U}}}\) for Good played according to the strategy. This means that there is a condition $p^*$ such that for every $\alpha<\lambda$, $p^*\leq p_\alpha$.  
\begin{claim}
    $p^*\in\bigcap j_{\mathcal{U}}''\mathcal{D}_p(\mathbb{P})$.
\end{claim}
\begin{proof}
    Let $D\subseteq \mathbb{P}$ be dense open above $p$. Let $C_0$ be the set of all $\alpha$ such that  $\mathbb{P}_\alpha=N_\alpha\cap \lambda\in \lambda$. Then by $(3)$ and by the continuity at points of cofinality $\geq\kappa$, $C_0$ is a club relative to cofinality $\geq \kappa$. Let $C_1$ be the collection of $\alpha$ such that $\mathbb{P}_\alpha=\alpha$. Then $C_1$ is still a club relative to cofinality $\geq\kappa$. There is a club $C_2$ such that for every $\alpha\in C_2$, $D\cap\mathbb{P}_\alpha$ is dense in $\mathbb{P}_\alpha$ (Take $C_2$ to be the club of closure points of the function mapping $\zeta<\lambda$ to an index $\zeta'\in D$ such that $\zeta\geq_{\mathbb{P}} \zeta'$).
  
Since \(\langle A_\alpha \mid \alpha \in S^{\lambda}_{\geq\kappa}\rangle\) is a diamond sequence, there is some $\alpha\in S^{\lambda}_{\geq\kappa}\cap C_0\cap C_1\cap C_2$ such that $D\cap \mathbb{P}_\alpha=D\cap \alpha =A_\alpha$. Therefore \(\alpha\) is an active stage, so by our construction, $q_\alpha\in j_{\mathcal{U}}(D\cap \mathbb{P}_\alpha)$. Then   $p^*\leq q_\alpha$, and therefore $p^*\in j_{\mathcal{U}}(D)$.
\end{proof}

\end{proof}
\begin{proposition}\label{Proposition in btween class}
    Let $\kappa$ be supercompact. Suppose that $\mathbb{P}$ is a forcing of cardinality $\lambda$ that is a regular cardinal which is not the successor of a singular of cofinality ${<}\kappa$. Then $\mathbb{P}$ is ${<}\kappa$-strategically closed implies $\mathbb{P}$ is ${<}\kappa$-diagonally distributive implies $\mathbb{P}$ is $\kappa$-distributive. 
\end{proposition}
\begin{proof}
    Suppose that $\mathbb{P}$ is $<\kappa$-strategically closed, and let us prove that it is ${<}\kappa$-diagonally distributive. Let $p\in\mathbb{P}$ and let $\vec{N}=\l N_i\mid i<\gamma\r$ be a $(\kappa,p,\mathbb{P})$-playful sequence. By elementarity, and since $H_\theta\models \mathbb{P}\text{ is }{<}\kappa\text{-strategically closed}$, there is $\vec{\sigma}=(\sigma_i\mid i<\kappa)\in N_0$ such that each $\sigma_i$ is a winning strategy for Good in the game of length $i$ with starting condition $p$. 
 
    Suppose that $\tau$ is a strategy for $G(\vec{N},\mathbb{P}/p)$. Let us show that $\tau$ is not winning for Bad. The point is to  play $\tau$ against $\sigma_{\gamma+1}$, and produce a winning run for Good. Note that since $\kappa+1\cup\{\vec{\sigma}\}\subseteq N_0$, $\sigma_{\gamma+1}\in N_0$. Start with $p_0=\sigma_{\gamma+1}(\l\r)$. Note that $p_0\in N_{0}$ (since $\sigma_{\gamma+1}\in N_{0}$), and let $q_0=\tau(\l p_0\r)\in N_{0}$ (by the rules of the game). Suppose that $\l p_i,q_i\mid i<\beta\r$ was defined, such that $p_i=\sigma_{\gamma+1}(\l p_j,q_j\mid j<i\r)$ and $q_i=\tau(\l p_j,q_j\mid j<i\r^\smallfrown p_i)$. Note that $(p_i,q_i)\in N_i$ for each $i<\beta$. Thus by condition $(ii)$, $\l p_i,q_i\mid i<\beta\r\in N_{\beta}$. Since $\sigma_{\gamma+1}\in N_{\beta}$, $\sigma_{\gamma+1}(\l p_i,q_i\mid i<\beta\r)=p_\beta\in N_{\beta}$ and also $q_\beta=\tau(\l p_i,q_i\mid i<\beta\r^\smallfrown p_\beta)\in N_{\beta}$. Since $\sigma_{\gamma+1}$ is a winning strategy for Good, there is a condition stronger than all those played before stage $\gamma$. Thus Good wins the run of length $\gamma$ against $\tau$, showing that $\tau$ is not winning for Bad.    

Next, suppose that $\mathbb{P}$ is ${<}\kappa$-diagonally distributive and let us show it is $\kappa$-distributive.  Let $\vec{D}=(D_\alpha\mid \alpha<\gamma)$ be a collection of $\gamma<\kappa$ many dense open subsets of $\mathbb{P}$. We would like to prove that $\bigcap_{i<\gamma}D_i$ is dense. Towards this, let $p\in\mathbb{P}$.  Let $\vec{N}$ be a $(\kappa,p,\mathbb{P})$-playful sequence such that $\vec{D},p\in N_0$. Note that by elementarity, for each $\alpha<\gamma$, $D_\alpha\cap N_\alpha$ is dense in $\mathbb{P}\cap N_\alpha$. 
Let us define a strategy for Bad in $G(\vec{N},\mathbb{P}/p)$ (which must then be non-winning).

The strategy is defined such that at stage $\alpha<\gamma$, if $(p_i,q_i\mid i<\alpha)$ is a partial run and $p_\alpha\in N_{\alpha}$ was the response of Good, then Bad chooses $q_\alpha\leq p_\alpha$ inside $D_\alpha\cap N_{\alpha}$. This is possible, since $p_\alpha\in N_\alpha\cap \mathbb{P}$ and $D_\alpha\cap N_\alpha$ is dense in $\mathbb{P}\cap N_{\alpha}$. By the assumption there is a winning run according to this strategy for Good which is below $p$. This produces a lower bound $p^*\in \bigcap_{\alpha<\gamma}D_\alpha$ such that $p^*\leq p$.

\end{proof}
\begin{proposition}\label{Proposition stationarity and diagonal distributivity}
    Let $\kappa<\lambda$ be regular cardinals such that $\delta^{<\kappa}<\lambda$ for every $\delta<\lambda$. If $\mathbb{P}$ is ${<}\kappa$-diagonally distributive, then $\mathbb{P}$ preserves stationary subsets of $S^\lambda_{<\kappa}$.
\end{proposition}
\begin{proof}
    Let $S$ be a stationary subset of $S^\lambda_{<\kappa}$, and suppose towards contradiction that $p\in\mathbb{P}$ and $\dot{C}$ is a $\mathbb{P}$-name for a club forced by $p$ to be disjoint from $\check{S}$. Let $\theta$ be sufficiently large and regular, and let $\vec{N}$ be a $\lambda$-sequence of elementary submodels such that:
    \begin{enumerate}
        \item [(a)] $\vec{N}$ is increasing and continuous.
        \item [(b)] $N_i\prec H_\theta$.
        \item [(c)] $|N_i|<\lambda$
        \item [(d)] $(N_i)^{<\kappa}\cup\{N_i\}\subseteq N_{i+1}$
        \item [(e)] $N_i\cap \lambda\in \lambda$.
        \item [(f)] $\kappa+1\cup\{\dot{C},p,\mathbb{P},\lambda\}\subseteq N_0$.
    \end{enumerate}
    To see that such a sequence exists, we define at successor steps $N_{i+1}\prec H_\theta$ such that $(a)-(f)$ is satisfied: $(d)$ is possible since $|N_i|^{<\kappa}<\lambda$. To ensure $(e)$ we iterate $\omega$-many times before defining $N_{i+1}$. At limits we must take unions by $(a)$. Note that $(a),(b),(c),(e)$ are maintained by unions, while $(d),(f)$ do not say anything about them.
    It follows that for every $\gamma<\kappa$ and a nowhere continuous increasing sequence $(\alpha_i\mid i<\gamma)$, $\vec{N}\restriction \{\alpha_i\mid i<\gamma\}$ is $(\kappa,p,\mathbb{P})$-playful: to see condition $(ii)$, for every $i$, the $i^{\text{th}}$ model here is $N_{\alpha_i}$, and $\alpha_i\geq\sup_{j<i}\alpha_j+1$. Hence $(\bigcup_{j<i}N_{\alpha_j})^{<\kappa}=(N_{\sup_{j<i}\alpha_j})^{<\kappa}\subseteq N_{\sup_{j<i}\alpha_j+1}\subseteq N_{\alpha_i}$ by $(d)$.

    Since $\vec{N}$ is continuous, and by $(e)$, $\{N_i\cap \lambda\mid i<\lambda\text{ limit}\}$ is a club and therefore there is $\rho<\lambda$ limit such that $\mu:=N_\rho\cap\lambda\in S$. Since $S\subseteq cof(<\kappa)$, $\cf(\mu)=\cf(\rho)=\gamma<\kappa$. Let $\{\alpha_i\mid i<\gamma\}$ be a nowhere continuous increasing and cofinal sequence in $\rho$. Then $\{N_{\alpha_i}\cap \lambda\mid i<\gamma\}$ is cofinal in $\mu$.

    Define a strategy in the game $G(\vec{N}\restriction \{\alpha_i\mid i<\gamma\},\mathbb{P}/p)$ by requiring that after Good plays $p_i$, Bad plays $q_i\in N_{\alpha_i}$ with $q_i\leq p_i$, forcing some ordinal $\xi_i$ into $\dot{C}$ above $N_{\sup_{j<i}\alpha_j}\cap\lambda$. Note that $p_i,\dot{C},N_{\sup_{j<i}\alpha_j}\in N_{\alpha_i}$. By elementarity, we may choose $\xi_i<N_{\alpha_i}\cap\lambda$. Since this is a non-winning strategy, there is a run $(p_\alpha,q_\alpha\mid \alpha<\gamma)$ according to this strategy which is winning for Good. In particular there is a lower bound $p^*\leq p_\alpha$ for all $\alpha$. Then $p^*\leq p$ and $p^*\Vdash \check{\mu}=\sup\{\check{\xi}_\alpha\mid \alpha<\gamma\}\in \check{S}\cap \dot{C}$, contradiction.
\end{proof}
This class is strictly in between the ${<}\kappa$-strategically closed forcings and $\kappa$-distributive forcings. For example, the forcing for shooting a club through a fat stationary subset of $\kappa$ is $\kappa$-distributive non-${<}\kappa$-strategically closed of size $\kappa$ (assuming $\kappa$ is a strong limit). Hence it is also ${<}\kappa$-diagonally distributive. 

In the next section we build a $\kappa$-distributive forcing which is not a projection and therefore not ${<}\kappa$-diagonally distributive.
\section{Preservation of Stationary sets}
In this section we construct examples of $\kappa$-distributive forcings of size $\lambda$ which are not projections of the $\lambda$-supercompact Prikry forcing. By Proposition \ref{Proposition stationarity and diagonal distributivity}, a natural example is the forcing for shooting a club through a stationary subset of $S^{\lambda}_{<\kappa}$. Before we do that, let us present a general result regarding projections of the Prikry forcing and the preservation of stationary sets:
\begin{theorem}\label{theorem: preservation of stationary}
    Suppose that $\mathcal{U}$ is a normal measure over $P_\kappa(\lambda)$. Suppose that $N$ is an intermediate model of $V^{Pr(\mathcal{U})}$ such that $N^{<\kappa}\subseteq V$ and in which $\lambda$ remains regular. Then $N$ preserves stationary subsets of $S^\lambda_{<\kappa}$. 
\end{theorem}
Suppose towards contradiction that $\mathbb{B}$ is a complete $\kappa$-distributive Boolean subalgebra of $\ro(Pr(\mathcal{U}))$, $S\subseteq S^\lambda_{<\kappa}$ is stationary in $\lambda$,  and $\dot{C}$ is a $\mathbb{B}$-name such that some condition $b\in\mathbb{B}$ forces $\dot{C}$ is a club in $\lambda$ which is disjoint from $\check{S}$ and $\check{\lambda}$ is a regular cardinal. We can find $p\in Pr(\mathcal{U})$ which forces $\check{b}$ in the induced generic for $\mathbb{B}$. Let us transfer this information to the $\omega^{\text{th}}$-iterated ultrapower. By Theorem \ref{Bukovsky-Dehornoy}, $\bigcap_{n<\omega} M_n$ is an $M_\omega$-generic extension by $j_\omega(Pr(\mathcal{U}))$, using a generic which includes $j_\omega(p)$.  Let $H$ be the induced $M_\omega$-generic for $j_{\omega}(\mathbb{B})$, so that $j_\omega(b)\in H$. By elementarity, in $M_\omega[H]$, $j_\omega(\lambda)$ is regular and  $j_\omega(\dot{C})_H$ is a club disjoint from $j_\omega(S)$. 
\begin{lemma}\label{Central Lemma}
     Let $\lambda\geq\kappa$ be regular, let $\mathcal{U}$ be a normal fine ultrafilter on $P_\kappa(\lambda)$, and let $\langle M_n,j_{m,n}\mid m\leq n\leq\omega\rangle$ be its iterated ultrapower.
     Suppose $\mathbb{Q}$ is a forcing and $H\in V$ is an $M_\omega$-generic filter on $j_{0,\omega}(\mathbb{Q})$ such that $j_{0,\omega}(\lambda)$ is regular in $M_\omega[H]$.
     Then there is $n<\omega$ such that $\sup\bigl(Hull^{M_\omega[H]}(j_{n,\omega}[M_n]\cup\{H\})\cap j_\omega(\lambda)\bigr)=\sup j_{n,\omega}''j_n(\lambda)$.
\end{lemma}
\begin{proof}[Proof of Lemma]
    Let $$M=Hull^{M_\omega[H]}(j_{0,\omega}[V]\cup\{H\}),$$ and let $\pi:M\to N$ be the transitive collapse and $$k=\pi^{-1},\text{ and }i=\pi\circ j_{0,\omega}:V\to N.$$ Denote by $\bar{M}=\pi[M_{\omega}\cap M]$ and $\bar{H}=\pi(H)$. Note that $N=\bar{M}[\bar{H}]$. Indeed, every $x\in N$ is $\pi(y)$ and $y=j_{0,\omega}(f)(H)$ for some $f\in V$. Hence $x=\pi(y)=i(f)(\bar{H})$. Note that $i(f)\in \bar{M}$ and therefore $x\in \bar{M}[\bar{H}]$. The other inclusion is similar. We also claim that $N$ is closed under $\kappa$-sequences from $V$. Indeed, if $(x_i\mid i<\kappa)\in N^\kappa\cap V$, choose $y_i$ such that $\pi(y_i)=x_i$.  Let $f_i\in V$ be such that $j_\omega(f_i)(H)=y_i$. Then $(j_\omega(f_i)\mid i<\kappa)=j_{\omega}((f_i\mid i<\kappa))\restriction \kappa$ is in $M$. Hence $(y_i\mid i<\kappa)=(j_\omega(f_i)(H)\mid i<\kappa)\in M$. It follows that $N$ is closed under $\kappa$-sequences from $V$.
    
    Since $j_{\omega}(\lambda)$ is regular in $M_{\omega}[H]$, also $M\models j_\omega(\lambda)$ is regular. Hence $i(\lambda)$ is regular in $N$. It follows that $i(\lambda)$ must have uncountable cofinality in $V$(!) Since $j_\omega(\lambda)$ has countable cofinality in $V$, $\nu:=\sup k[i(\lambda)]<k(i(\lambda))=j_{\omega}(\lambda)$. Let us argue that $$\nu=j_\omega(\lambda)\cap Hull^{M_\omega[H]}(k[N]\cup\nu).$$ Clearly $\nu\leq j_\omega(\lambda)\cap Hull^{M_\omega[H]}(k[N]\cup\nu)$. Also, if $\rho\in Hull^{M_\omega[H]}(k[N]\cup\nu)\cap j_\omega(\lambda)$, then there is $f\in V$ and $\vec{\xi}\in [\nu]^{<\omega}$ such that $j_{\omega}(f)(H,\vec{\xi})=\rho$. We may assume that $\text{rng}(f)=\lambda$. Consider $i(f)\in N$, and consider $\bar{\rho}<i(\lambda)$ such that $\vec{\xi}\leq k(\bar\rho)$. Since $i(\lambda)$ is regular, $\mu=\sup_{\vec{\zeta}\in [\bar\rho]^{<\omega}}i(f)(\bar{H},\vec{\zeta})<i(\lambda)$. Then $\nu\geq k(\mu)=\sup_{\vec{\zeta}\in [k(\bar\rho)]^{<\omega}} j_{\omega}(f)(H,\vec{\zeta})\geq\rho$. 
    Let $\xi_n=\sup j_{n,\omega}''j_n(\lambda)$. Then $\{\xi_n\mid n<\omega\}$ is cofinal in $j_\omega(\lambda)$ and therefore there is a least $n$ such that $\nu\leq\xi_n$. Every $\rho<\xi_n=\sup j_{n,\omega}''j_n(\lambda)$ is bounded by some $j_{n,\omega}(\bar\rho)$ where $\bar\rho<j_n(\lambda)$. But then there is $g$ with $\text{rng}(g)=\lambda$ such that $j_n(g)(j_{n,\omega}^{-1}\xi_0,...,\xi_{n-1})=\bar\rho$ and $j_\omega(g)(\xi_0,\ldots,\xi_{n-1})=j_{n,\omega}(\bar\rho)$. Hence $j_{n,\omega}(\bar\rho)\in Hull^{M_\omega[H]}(j_{0,\omega}[V]\cup \nu)\cap j_\omega(\lambda)=\nu$, so $\rho<\nu$.
    
\end{proof}
\textit{End of the proof of Theorem \ref{theorem: preservation of stationary}.} Let $n<\omega$ be as in Lemma \ref{Central Lemma}. We note that $Hull^{M_\omega[H]}(j_{n,\omega}[M_n]\cup\{H\})$ is closed under $j_n(\kappa)$-sequences from $M_n$. Indeed, suppose that $\langle x_\alpha\mid \alpha<j_n(\kappa)\rangle\subseteq Hull^{M_\omega[H]}(j_{n,\omega}[M_n]\cup\{H\})$ and $\langle x_\alpha\mid \alpha<j_n(\kappa)\rangle\in M_n$. Pick $\langle y_\alpha\mid \alpha<j_n(\kappa)\rangle\in M_n$ such that $j_{n,\omega}(y_\alpha)(H)=x_\alpha$. Since the critical point of $j_{n,\omega}$ is $j_n(\kappa)$, we conclude that $$\langle j_{n,\omega}(y_\alpha)\mid \alpha<j_n(\kappa)\rangle=j_{n,\omega}(\langle y_\alpha\mid \alpha<j_{n}(\kappa)\rangle)\restriction j_{n}(\kappa)\in Hull^{M_\omega[H]}(j_{n,\omega}[M_n]\cup\{H\}),$$ and also
$$\langle x_\alpha\mid \alpha<j_n(\kappa)\rangle=\l j_{n,\omega}(y_\alpha)(H)\mid \alpha<j_n(\kappa)\rangle\in Hull^{M_\omega[H]}(j_{n,\omega}[M_n]\cup\{H\}).$$

We let $C:=j_{\omega}(\dot{C})_H\in Hull^{M_\omega[H]}(j_{n,\omega}[M_n]\cup\{H\})$ be a club disjoint from $j_{\omega}(S)$. Since $\sup\bigl(Hull^{M_\omega[H]}(j_{n,\omega}[M_n]\cup\{H\})\cap j_\omega(\lambda)\bigr)=\sup j_{n,\omega}''j_n(\lambda)$, we conclude that $C\cap \sup j_{n,\omega}''j_n(\lambda)$ is a club in $\sup j_{n,\omega}''j_n(\lambda)$. Hence $C\cap j_{n,\omega}''j_n(\lambda)$ is a ${<}j_n(\kappa)$-closed unbounded subset of $\sup j_{n,\omega}''j_n(\lambda)$. By closure under $<j_n(\kappa)$-sequences, $C^*:=j_{n,\omega}^{-1}[C\cap j_{n,\omega}''j_n(\lambda)]$ is a $<j_n(\kappa)$-closed subset of $j_n(\lambda)$ which is disjoint from $j_n(S)$. This contradicts that $j_{n}(S)$ is a stationary subset of $(S^{j_n(\lambda)}_{<j_n(\kappa)})^{M_n}$.\qed

Next we would like to force a situation where $\kappa$ is supercompact and there is a forcing which yields the situation of Theorem \ref{theorem: preservation of stationary}. As we said previously, we will consider the forcing for shooting a club through a stationary set:
\begin{definition}
    Let $\lambda$ be a regular cardinal and let $S\subseteq \lambda$ be stationary. The forcing $Q(S)$ consists of conditions $a\subseteq S$ which are closed, bounded in $\lambda$ and $\max(a)$ exists. The order on $Q(S)$ is end extension. 
\end{definition}
\begin{definition}
    Let $\theta\leq\lambda$ be regular cardinals. A stationary set $S\subseteq \lambda$ is called $\theta$-fat if for any club $C\subseteq \lambda$, and every $\tau<\theta$, $S\cap C$ contains a closed set of order-type $\tau+1$. It is called \textit{fat} if it is $\lambda$-fat.
\end{definition}    
   The following was essentially proven by Abraham and Shelah \cite{AbrahamShelah1994}.
\begin{proposition}\label{abraham-shelah}
    Let $\theta\leq\lambda$ be regular and suppose that $\lambda$ is either inaccessible or $\lambda=\kappa^+$ with $\kappa^{<\kappa}=\kappa$. If  $S\subseteq\lambda$ is $\theta$-fat, then the forcing shooting a club through $S$ is $\theta$-distributive.
\end{proposition}
\begin{proof}
    Let $\vec{D}=\l D_i\mid i<\rho\r$ for $\rho<\theta$. We let $\langle M_i\mid i<\lambda\rangle$ be a sequence of elementary models which is increasing and continuous, $\rho\cup\{\vec{D},\theta,\rho,\lambda\}\subseteq M_0$, $M_{i+1}$ contains $<\kappa$-sequences from $M_i$, and $M_i\cap\lambda\in\lambda$. Such a sequence exists by the cardinal arithmetic assumptions on $\lambda$. Let $C$ be the club  $\{M_i\cap\lambda\mid i<\lambda\}$ and find a closed set $A\subseteq S\cap C$ of order type $\rho+1$. Then we can recursively meet the dense sets $D_i$ with a condition $p_i$ so that $(\sup(p_i),\sup(p_{i+1}))\cap A\neq\emptyset$. At the end we can take the union of the $p_i$ and add the supremum which is guaranteed to be in $A$ (and so in $S$) forming a legitimate condition in $Q(S)$.
\end{proof}
\begin{corollary}
    Suppose that the supercompactness of $\kappa$ is indestructible under $\Add(\kappa^+,1)$. Let $G\subseteq \Add(\kappa^+,1)$ be $V$-generic. Then in $V[G]$ there is a fat stationary set $S\subseteq\kappa^+$ such that $S^{\kappa^+}_{<\kappa}\setminus S$ is stationary. In particular, $Q(S)$ cannot be a projection of $Pr(\mathcal{U})$ for any normal ultrafilter $\mathcal{U}$ over $P_\kappa(\kappa^+)$.
\end{corollary}
\begin{proof}
    Let $A_G\subseteq\kappa^+$ be the generic $\Add(\kappa^+,1)$ set. Let us show that $S^{\kappa^+}_{<\kappa}\setminus A_G$ is stationary. Let $\dot{C}$ be an $\Add(\kappa^+,1)$-name for a club and suppose towards a contradiction that some $p$ forces $(S^{\kappa^+}_{<\kappa}\setminus A_G)\cap \dot{C}=\emptyset$. We construct an $\omega$-sequence $$p=p_0\geq p_1\geq p_2\geq ...$$ 
    such that $p_{n+1}$ forces some ordinal $\xi_n\in (\sup(\text{supp}(p_n)),\kappa^+)\cap \dot{C}$. Then letting $\xi^*=\sup_{n<\omega}\xi_n$ and $p'=\bigcup_{n<\omega}p_n$, we let $p^*=p'\cup \{(\xi^*,0)\}$. Then $cf(\xi^*)=\omega$ and $p^*\Vdash \xi^*\in S^{\kappa^+}_{<\kappa}\setminus A_{\dot{G}}$. Also $p^*\Vdash\xi^*\in \dot{C}$, and $p^*\Vdash \xi^*\in (S^{\kappa^+}_{<\kappa}\setminus A_{\dot{G}})\cap \dot{C}$. On the other hand $p^*\leq p$, contradiction.

    Finally, we note that $A_G$ is $\kappa^+$-fat. This is a well-known fact whose short proof we include here for the sake of completeness. Indeed, if $\dot{C}$ is a $\Add(\kappa^+,1)$-name forced by a condition $p$ to be a club, and $\tau<\kappa^+$, we proceed with a density argument. Let $p'\leq p$, and construct a decreasing sequence of conditions $(p_i\mid i\leq\omega\cdot\tau)$ such that:
    \begin{enumerate}
        \item for every $i<\omega\cdot\tau$, there is some ordinal $\xi_i\in (\sup(\text{supp}(p_i)),\kappa^+)$ such that $p_{i+1}\Vdash \xi_i\in \dot{C}$.
        \item For a limit ordinal $i\leq\omega\cdot\tau$, $p_i=(\bigcup_{j<i}p_j)\cup\{(\sup_{j<i}\xi_j,1)\}$.
    \end{enumerate}
    Such a sequence can easily be constructed recursively using the $\kappa^+$-closure of $\Add(\kappa^+,1)$. It follows that $p_{\omega\cdot\tau}\Vdash\{\sup_{j<\omega\cdot i}\xi_j\mid  i\leq\tau\}\subseteq \dot{C}\cap A_{\dot{G}}$. 
\end{proof}
In \cite{Krueger}, Krueger calls sets $S\subseteq \kappa^+$ which are fat and $cf(<\kappa)\setminus S$ being stationary, non-trivial fat stationary sets. Krueger proves that in $L$ (or if $\square_\kappa$ holds) there are non-trivial fat stationary subsets of $\kappa^+$. On the other hand, he shows that if every stationary subset $S\subseteq S^{\kappa^+}_{<\kappa}$ reflects on a club relative to $S^{\kappa^+}_\kappa$, then there are no non-trivial fat stationary subsets of $\kappa^+$. Assuming the existence of a weakly compact cardinal, Magidor \cite{Magidor1982} constructed a model where this type of stationary reflection holds at the successor $\kappa^+$ of a regular cardinal $\kappa$. This model is obtained by collapsing a weakly compact $\lambda>\kappa$ using the usual L\'{e}vy collapse $Col(\kappa,<\lambda)$, followed by a $<\kappa$-support iteration of length $\kappa^{++}$ of shooting $\kappa$-clubs through fat stationary subsets of $\kappa^+$. This forcing is trivially ${<}\kappa$-directed closed. The important point, which was proven by Magidor, is that $\kappa^+$ is preserved by the iteration. Hence, if we start with a Laver-indestructible supercompact cardinal $\kappa$ \cite{Laver1978} with a weakly compact cardinal $\lambda>\kappa$, we obtain the following corollary:  
\begin{corollary}
   Assume that there is a supercompact cardinal $\kappa$ and a weakly compact cardinal $\lambda>\kappa$. Then it is consistent that  $\kappa$ is supercompact and that every fat stationary subset of $\kappa^+$ is trivial.
\end{corollary}
It is still possible that there is a $\kappa^+$-distributive forcing of size $\kappa^+$ which is not of the form $Q(S)$ that cannot be a projection of $Pr(\mathcal{U})$, hence the following question:
\begin{question}
    Is it consistent that every $\kappa^+$-distributive forcing of size $\kappa^+$ is a projection of the supercompact Prikry forcing?
\end{question}
Relaxing the $\kappa^+$-distributivity, our next goal is to show that in $ZFC$ we can always find a $\kappa$-distributive forcing of size $\kappa^+$ which is not a projection of $Pr(\mathcal{U})$ for any normal ultrafilter $\mathcal{U}$. Again, we will find suitable stationary sets. This time we will use $\kappa$-fat stationary sets rather than fully fat stationary sets. For these, Krueger has the following $ZFC$ result \cite[Theorem 2]{Krueger}:
\begin{theorem}[Krueger]
Let  $\mu=\kappa^+$ where $\kappa$ is a regular uncountable cardinal and let $A\subseteq \mu$ be a $\kappa$-fat stationary set. Then every stationary subset $B\subseteq A$ contains a stationary set $B_0$ such that $A \setminus B_0$ is $\kappa$-fat. 
\end{theorem}
In particular if $\kappa$ is supercompact, then taking $A=S^{\kappa^+}_{<\kappa}$, there is a stationary subset $S^*\subseteq S^{\kappa^+}_{<\kappa}$ such that $S^{\kappa^+}_{<\kappa}\setminus S^*$ is $\kappa$-fat. Then, by  Proposition \ref{abraham-shelah}, the forcing $Q(S^{\kappa^+}_{<\kappa}\setminus S^*)$ which shoots a club through $S^{\kappa^+}_{<\kappa}\setminus S^*$ is $\kappa$-distributive. Note that since the forcing kills the stationarity of $S^{\kappa^+}_\kappa$ (or simply because it adds a club of order type $\kappa$ to $\kappa^+$), it collapses $\kappa^+$. 

\begin{remark}
    Note that since the forcing is $\kappa$-distributive, a generic filter added by the above forcing is closed under ${<}\kappa$-sequences from the ground model.
\end{remark}
\begin{theorem}\label{theorem: kappa fat no projection}
    Let $\mathcal{U}$ be a normal ultrafilter on $P_\kappa(\kappa^+)$, and let $S\subseteq \kappa^+$ be a $\kappa$-fat stationary subset of $\kappa^+$ such that $S^{\kappa^+}_{<\kappa}\setminus S$ is also stationary. Then there is no projection of $Pr(\mathcal{U})$ onto $Q(S)$.
\end{theorem}
\begin{proof}
    By $\kappa$-fatness of $S$, Proposition \ref{abraham-shelah} entails that $Q(S)$ is $\kappa$-distributive. Hence by Theorem \ref{Projection criteria}, to rule out a projection, it suffices to show that for every $n<\omega$, $\bigcap j_{\mathcal{U}^n}''\mathcal{D}(Q(S))=\emptyset$. Suppose otherwise that $\bigcap j_{\mathcal{U}^n}''\mathcal{D}\neq\emptyset$. For now, let us prove a series of propositions for an arbitrary elementary embedding $j:V\to M$ with critical point $\kappa$ such that $M$ is closed under $\kappa$-sequences such that $\bigcap j''\mathcal{D}(Q(S))\neq\emptyset$. 
\begin{center}
Suppose that $q=q(j)\in \bigcap j''\mathcal{D}(Q(S))$ has minimal $\max(q)$.
\end{center}
Then for each $\eta<\max(q)$, there is a dense open $D\subseteq Q(S)$ such that $q\cap \eta\notin j(D)$. For the next proposition, it is essential that $S$ is $\kappa$-fat:
\begin{proposition}\label{prop: cutting initial segment}
     For every $\eta<\max(q)$, $q\setminus\eta\in \bigcap j''\mathcal{D}(Q(S))$.
\end{proposition}
 \begin{proof}
     Otherwise, there are $D_1$, $\eta_1$ such that $q\setminus \eta_1\notin j(D_1)$. By minimality, there is $D_2$ such that $q\cap(\eta_1+1)\notin j(D_2)$. Define $D^*$ as the set of all $p$ such that for some  $\eta<\max(p)$:
     \begin{enumerate}
         \item [(I)] $p\cap(\eta+1)\in D_2$.
         \item [(II)] for every $\alpha\leq\eta$, $p\setminus\alpha\in D_1$.
     \end{enumerate}
     This is easily seen to be open. To see that it is dense, let $p$ be any condition. Extend $p$ to some $p_1\in D_2$. 
     We fix a sequence $(M_\alpha\mid \alpha<\kappa^+)$ and $\theta$ large enough so that:
     \begin{enumerate}
         \item for each $\alpha<\kappa^+$, $M_\alpha\prec H(\theta)$.
         \item $|M_\alpha|=\kappa$.
         \item The sequence $M_\alpha$ is increasing and continuous.
         \item $\{p_1,S,D_1\}\cup\kappa+1\cup \sup(p_1)\subseteq M_0$.
         \item $\gamma_\alpha:=M_\alpha\cap\kappa^+\in \kappa^+$.
         \item $M_\alpha\in M_{\alpha+1}$ and $M_\alpha^{<\kappa}\subseteq M_{\alpha+1}$
     \end{enumerate} Let $C=\{\gamma_\alpha\mid \alpha<\kappa^+\}$. Let $\text{otp}(p_1)=\rho<\kappa$. By fatness of $S$, $S\cap C$ contains a closed set $A$ of order-type $\rho+1$. Enumerate increasingly $A=\{\gamma_{\alpha_i}\mid i\leq\rho\}$, and $p_1=\{\delta_i\mid i\leq\rho\}$. Now we define recursively a decreasing sequence $r_i$ for $i\leq\rho$ such that:
     \begin{enumerate}
         \item [(a)] $r_i\in M_{\alpha_{i+1}}$.
         \item [(b)] $\max(r_i)\in [\gamma_{\alpha_i},\gamma_{\alpha_{i+1}})$.
         \item [(c)] $r_{i+1}\setminus (\delta_i+1)\in D_1$.
     \end{enumerate}
     
     Let $p_1\geq r_0\in D_1\cap M_{\alpha_1}$, with $\max(r_0)\geq\gamma_{\alpha_0}$. Such an $r_0$ exists, since $M_{\alpha_1}\models ``D_1$ is dense" so $p_1$ has an extension $r_0'\in D_1\cap M_{\alpha_1}$ and also $M_{\alpha_1}\models S\setminus\gamma_{\alpha_0}\neq\emptyset$, hence if necessary, we can add a point from $S\cap M_{\alpha_1}$ to $r'_0$ to form $r_0\in M_{\alpha_1}\cap D_1$ such that $\max(r_0)\in [\gamma_{\alpha_0},\gamma_{\alpha_1})$. 
     
     At successor steps, suppose that $r_i\in M_{\alpha_i+1}$ was defined. We let $r_{i+1}'\leq r_i$ be an element of $M_{\alpha_{i+2}}$ such that $\max(r'_{i+1})\geq\gamma_{\alpha_{i+1}}$. Then, let $r_{i+1}''\leq r'_{i+1}\setminus (\delta_i+1)$ in $M_{\alpha_{i+2}}\cap D_1$ and let $$r_{i+1}=r_{i+1}''\cup (p_1\cap (\delta_i+1))\in M_{\alpha_{i+2}}.$$ Note that
     since $r_i\leq p_1$, and since $r_{i+1}''\leq r_i\setminus (\delta_i+1)$, we have $r_{i+1}\leq r_i$. Clearly, $r_{i+1}$ satisfies $(a),(b),(c)$.
     
     At limit steps $\mu$, we have constructed a decreasing sequence $r_i\in M_{\alpha_{i+1}}\subseteq M_{\alpha_\mu}$. By condition $(6)$, the sequence $\{r_i\mid i<\mu\}\in M_{\alpha_{\mu+1}}$ and therefore,  $\bigcup_{i<\mu}r_i\in M_{\alpha_{\mu+1}}$. By condition $(b)$, we have $$\sup(\bigcup_{i<\mu}r_i)=\gamma_{\alpha_\mu}\in A\cap M_{\alpha_\mu}\subseteq S\cap M_{\alpha_{\mu+1}},$$ so we can form the condition  $r'_{\mu}=\bigcup_{i<\mu}r_i\cup\{\gamma_{\alpha_\mu}\}$. As in the successor step, we extend $r'_{\mu}\setminus (\delta_\mu+1)$ to a condition $r''_{\mu}\in D_1\cap M_{\alpha_{\mu+1}}$, and set $r_\mu=r''_\mu\cup (p_1\cap (\delta_\mu+1))$. This concludes the definition of the sequence of the $r_i$'s.
     
     Let $p_2:=r_{\rho}$. Then $p_2\leq p_1\leq p$ and we claim that $p_2\in D^*$ as witnessed by $\eta=\max(p_1)$. For each $\alpha\leq\max(p_1)$, there is $i<\rho$ such that $p_1\setminus \alpha=p_1\setminus (\delta_i+1)$. The condition $p_2\setminus \alpha$ is then stronger than $r_{i+1}\setminus \alpha=r_{i+1}\setminus (\delta_i+1)\in D_1$ (by condition $(c)$). Hence also $p_2\setminus\alpha\in D_1$.
     
      So far we have shown that $D^*$ is dense open, and therefore $q\in j(D^*)$. Let $\eta<\max(q)$ witness this. If $\eta\leq\eta_1$, then since $q\cap(\eta_1+1)\notin j(D_2)$, and $D_2$ is open, $q\cap(\eta+1)\notin j(D_2)$ contradicting (I). If $\eta>\eta_1$, then by condition (II), $q\setminus\eta_1\in j(D_1)$, contradicting the choice of $D_1$. 
 \end{proof}
 \begin{corollary}
     $q$ is unbounded in $\max(q)$. In particular $q\cap \max(q)$ is a club in $j(S)\cap \max(q)$.
 \end{corollary}
 \begin{proof}
     Let $D_2$ be the dense open set saying that there are at least three elements. Since $q\in j(D_2)$, if $q\cap\max(q)$ were bounded, it would have a maximum $\beta<\max(q)$. Then $q\setminus\beta$ has exactly two elements, but $q\setminus\beta\in j(D_2)$, a contradiction.
 \end{proof}
 For the next proposition we only need that $S$ is unbounded in $\kappa^+$ and $S^{\kappa^+}_{<\kappa}\setminus S$ is stationary in $\kappa^+$:
 \begin{proposition}\label{max properties}
 $\max(q)>\sup j''\kappa^+$. Moreover,  whenever $i:V\to N$ and $k:N\to M$ are elementary, such that $k$ is definable in $N$, $j=k\circ i$, $q\in N$, and $\text{crit}(k)\geq i(\kappa)$, we can conclude that $\max(q)\neq \sup k''i(\kappa^+)$.
 \end{proposition}
 \begin{proof}
 For the first part, for each $\alpha<\kappa^+$ let $D_\alpha$ be the dense open set saying that there are elements above $\alpha$. This is dense since $S$ is unbounded. Since $q\in j(D_\alpha)$, $\max(q)>j(\alpha)$ and $\max(q)\geq \sup j''\kappa^+$. The strict inequality follows from the second part by taking $i=id$ and $j=k$.
     Otherwise, we note that $k''i(\kappa^+)$ is an ${<}i(\kappa)$-club in $\sup k''i(\kappa^+)=\max(q)$ and so is $q\cap \max(q)$. Let $C^*=k^{-1}[k''i(\kappa^+)\cap q]$. Then:
     \begin{enumerate}
        \item $C^*\in N$ (since $k$ is definable in $N$).
         \item $C^*\subseteq i(S)$.
         \item $C^*$ is a ${<}i(\kappa)$-club.
     \end{enumerate}
     This is a contradiction in $N$, since by elementarity of $i$, $S^{i(\kappa^+)}_{{<}i(\kappa)}\setminus i(S)$ is stationary in $i(\kappa^+)$ in $N$.
 \end{proof}
From now on we assume that $j=j_n=j_{\mathcal{U}^n}$. By normality of $\mathcal{U}$, we have that
$$[id]_{\mathcal{U}^n}=(j_n''\kappa^+,j_{1,n}''j_1(\kappa^+),...,j_{n-1,n}''j_{n-1}(\kappa^+)).$$
Moreover, $$j_n''\kappa^+\prec j_{1,n}''j_1(\kappa^+)\prec...\prec j_{n-1,n}''j_{n-1}(\kappa^+).$$
By Proposition \ref{max properties},  for every $0\leq i\leq n-1$, $\max(q)\neq \sup j_{i,n}''j_{i}(\kappa^+)$. Let $i^*<n$ be such that $\max(q)\in (\sup j_{i^*,n}''j_{i^*}(\kappa^+),\sup j_{i^*+1,n}''j_{i^*+1}(\kappa^+))$. Denote by $\delta^*=\sup j_{i^*,n}''j_{i^*}(\kappa^+)$.
Since $\delta^*<\max(q)$, by Proposition \ref{prop: cutting initial segment} we have: \begin{corollary}
    \label{definition of qstar}
     $q^*(j):=q\setminus (\delta^*+1)\in \bigcap j''D(Q(S)).$
\end{corollary} Let us show that we can make some more modifications below $q^*(j)$:
 \begin{proposition}
    Let $S$ be ${<}\kappa$-fat. For every $p\in Q(S)$, $j(p)\cup q^*(j)\in \bigcap j''D(Q(S))$.
\end{proposition}
\begin{proof}
    The proof is essentially homogeneity, but let us add some details. Fix $D_1$, $p_1$, and let us show that $j(p_1)\cup q^*(j)\in j(D_1)$. Define $D^*$ as the set of all $p$ such that $p_1\cup (p\setminus(\max(p_1)+1))\in D_1$. Suppose we extend $p\in D^*$ to some $p'$, then $p_1\cup (p\setminus(\max(p_1)+1))\sqsubseteq p_1\cup (p'\setminus(\max(p_1)+1))$, so $p'\in D^*$, hence $D^*$ is open. To see $D^*$ is dense, take any $p$, extend it to $p'$ with $\max(p')>\max(p_1)+1$. Find $q\in D_1$ such that $q\sqsupseteq (p'\setminus (\max(p_1)+1))\cup p_1$. Set $p''=q\cup (p'\cap(\max(p_1)+1))$. Then $p''\leq p'\leq p$, and $(p''\setminus(\max(p_1)+1))\cup p_1=q\in D_1$. 
     
      Hence $D^*$ is dense open, and by (\ref{definition of qstar}), $q^*(j)\in j(D^*)$. This means that  $j(p_1)\cup(q^*(j)\setminus(\max(j(p_1))+1))\in j(D_1)$. But $\max(j(p_1))+1<\delta^*$. Hence $j(p_1)\cup q^*(j)\in j(D_1)$. 
\end{proof}

Let
 $f:[P_\kappa(\kappa^+)]^n\to Q(S)$ represent $q^*$ in $M_{\mathcal{U}^n}$, i.e. $j_n(f)([id]_{\mathcal{U}^n})=q^*$. 
\begin{corollary}\label{cor: entering dense sets1}
   For every dense open $D\subseteq Q(S)$, $$E_D:=\{\vec{x}\in [P_\kappa(\kappa^+)]^n\mid \forall \vec{y}\prec \vec{x}_0, \ \bigcup_{k\leq \frac{|\vec{y}|}{n}} f(\vec{y}\restriction[kn,(k+1)n))\cup f(\vec{x})\in D\}\in \mathcal{U}^n.$$
\end{corollary}
\begin{proof}
    Let us show that $[id]_{\mathcal{U}^n}\in j_n(E_D)$. Fix any $\vec{y}=(y_1,...,y_m)\prec j_n''\kappa^+=([id]_{\mathcal{U}^n})_0$. Then for each $1\leq i\leq m$, there is $z_i$ such that $j_n(z_i)=j_n''z_i=y_i$, and $j_{n}(\vec{z})=\vec{y}$. This means that for each $k\leq \frac{|\vec{y}|}{n}$, $$j_n(f)(\vec{y}\restriction [kn,(k+1)n))=j_n(f)(j_n(\vec{z})\restriction [kn,(k+1)n))=j_n(f(\vec{z}\restriction [kn,(k+1)n))),$$ and $$\bigcup_{k\leq \frac{|\vec{y}|}{n}} j_n(f)(\vec{y}\restriction [kn,(k+1)n))=j_n(\bigcup_{k\leq \frac{|\vec{y}|}{n}} f(\vec{z}\restriction [kn,(k+1)n))).$$ Note that $p:=\bigcup_{k\leq \frac{|\vec{y}|}{n}} f(\vec{z}\restriction [kn,(k+1)n))$ is a closed subset of $S$ and therefore $p\in Q(S)$. By the above computation and the previous lemma, $$\bigcup_{k\leq \frac{|\vec{y}|}{n}}j_n(f)(\vec{y}\restriction [kn,(k+1)n))\cup j_n(f)([id]_{\mathcal{U}^n})=j_n(p)\cup q^*\in j_n(D),$$
    which by definition means $[id]_{\mathcal{U}^n}\in j_n(E_D)$.
\end{proof}
\begin{claim}
There is an increasing function $g:\kappa^+\to\kappa^+$ such that $j_n(g)(\delta^*)\geq \max(q^*)$. 
\end{claim}
\begin{proof}
Since $\max(q^*)<\sup(j_{i^*+1,n}''j_{i^*+1}(\kappa^+))$, there is $\gamma<j_{i^*+1}(\kappa^+)$ such that $\max(q^*)\leq j_{i^*+1,n}(\gamma)$. Work in $M_{i^*+1}$. Then there is $g':[P_\kappa(\kappa^+)]^{i^*+1}\to\kappa^+$ such that $\gamma=j_{i^*+1}(g')([id]_{\mathcal{U}^{i^*+1}})$ and recall that $[id]_{\mathcal{U}^{i^*+1}}$ is $\prec$-increasing with a top element  $j_{i^*,i^*+1}''j_{i^*}(\kappa^+)$.
    For each $\alpha<\kappa^+$, set $g(\alpha)=\sup\{g'(\vec{y})\mid \vec{y}\in [P_\kappa(\alpha+1)]^{i^*+1}\}$. Since $\alpha^{<\kappa}<\kappa^+$, $g(\alpha)<\kappa^+$. Back in $M_{i^*+1}$ we have   $$j_{i^*+1}(g)(\sup(j_{i^*,i^*+1}''j_{i^*}(\kappa^+)))\geq j_{i^*+1}(g')([id]_{\mathcal{U}^{i^*+1}})=\gamma.$$ Applying $j_{i^*+1,n}$ to this inequality, we get
    $j_{n}(g)(j_{i^*+1,n}(\sup j''_{i^*,i^*+1}j_{i^*}(\kappa^+)))\geq j_{i^*+1,n}(\gamma)\geq\max(q^*)$.
    Finally, we note that $j_{i^*+1,n}(\sup j''_{i^*,i^*+1}j_{i^*}(\kappa^+))=\delta^*$: indeed, in $M_{i^*+1}$, the ordinal $\sup j''_{i^*,i^*+1}j_{i^*}(\kappa^+)$ has cofinality $j_{i^*}(\kappa^+)<j_{i^*+1}(\kappa)=\crit(j_{i^*+1,n})$. Hence 
    $$j_{i^*+1,n}(\sup j''_{i^*,i^*+1}j_{i^*}(\kappa^+))=\sup(j_{i^*+1,n}[j_{i^*,i^*+1}[j_{i^*}(\kappa^+)]])=\sup j_{i^*,n}[j_{i^*}(\kappa^+)]=\delta^*.$$
    
\end{proof}
\begin{claim}
    $E=\{\vec{x}\in [P_\kappa(\kappa^+)]^n\mid \forall y\prec x_0,\ g(\sup y)<\sup x_0\}\in \mathcal{U}^n.$
\end{claim}
\begin{proof}
     Indeed if $y\prec j_n''\kappa^+$, there is $z$ such that $j_n(z)=y$ and $$j_n(g)(\sup y)=j_n(g(\sup z))<\sup(j_n''\kappa^+).$$
\end{proof}

By forcing below a measure one set, any Prikry generic for $\mathcal{U}$ yields a Prikry generic for $\mathcal{U}^n$, $(\vec{x}(m))_{m<\omega}$, where $\vec{x}(m)=(\vec{x}(m)_0,...,\vec{x}(m)_{n-1})$. The purpose of the measure one set is that the following inequalities are satisfied for each $m<\omega$:
 $$\sup(\vec{x}(m)_{i^*})<^{(1)}\min(f(\vec{x}(m)))\leq \max(f(\vec{x}(m)))\leq^{(2)} g(\sup(\vec{x}(m)_{i^*}))<^{(3)}\sup(\vec{x}(m+1)_0).$$
 $(1)$ follows since $\min(q^*)>\delta^*$.\\
 $(2)$ follows since $j_n(g)(\delta^*)\geq\max(q^*)$.\\
 $(3)$ follows from the previous claim.

Let $q_m=f(\vec{x}(m))$. The above shows in particular that the $q_m$'s form blocks. Let $C:=\bigcup_{m<\omega}q_m$.
\begin{proposition}
    $C$ is $V$-generic for $Q(S)$.
\end{proposition}
\begin{proof}
    Let $D$ be dense open. Then by the corollary \ref{cor: entering dense sets1}, there is $m$ such that $\vec{x}(m)\in E_D$. This means that for every $\vec{y}\prec \vec{x}(m)_0$, $\bigcup_{k\leq |\vec{y}|/n}f(\vec{y}\restriction [kn,(k+1)n))\cup f(\vec{x}(m))\in D$. In particular, since $\vec{x}(0)^\smallfrown\vec{x}(1)^\smallfrown...^\smallfrown \vec{x}(m-1)\prec \vec{x}(m)$, we get that
    $C\cap (\max(q_m)+1)=\bigcup_{i\leq m}q_i=\bigcup_{i\leq m}f(\vec{x}(i))\in D$. It follows that $C\cap (\max(q_m)+1)\in D\cap G(C)$.
\end{proof}

\textit{End of the Proof of Theorem \ref{theorem: kappa fat no projection}.} To reach a contradiction, we will show that the cofinality of $\kappa^+$ in $V[C]$ is $\omega$ (contradicting the $\kappa$-distributivity of the forcing). Let $\gamma_0=\sup(\vec{x}(0)_{i^*})$. Then $\delta_0=g(\gamma_0)$ and let $\gamma_1=\min(C\setminus \delta_0+1)$, let $\delta_1=g(\gamma_1)$. In general, if $\delta_m$ is defined, we let $\gamma_{m+1}=\min(C\setminus \delta_m+1)$ and $\delta_{m+1}=g(\gamma_{m+1})$. We claim that $\gamma_m$ (which is clearly defined in $V[C]$) is unbounded in $\kappa^+$:

It suffices to show that for every $m$, $\gamma_m\geq \sup(\vec{x}(m)_{i^*})$. By induction. For $m=0$, by definition $\gamma_0=\sup(\vec{x}(0)_{i^*})$. Suppose this is true for $m$. Then $\gamma_m\geq\sup(\vec{x}(m)_{i^*})$. Then $g(\gamma_m)=\delta_m\geq g(\sup(\vec{x}(m)_{i^*}))\geq\max(f(\vec{x}(m)))=\max(q_m)$ (by monotonicity and the choice of $g$). Finally, $$\gamma_{m+1}=\min(C\setminus\delta_m+1)\geq \min(C\setminus \max(q_m)+1)=\min(q_{m+1})\geq\sup(\vec{x}(m+1)_{i^*}),$$
as wanted. \end{proof}
\begin{question}
    Does every projection of the $\kappa^+$-supercompact Prikry forcing in which $cf((\kappa^+)^V)\geq\kappa$ preserve stationary subsets of $S^{\kappa^+}_{<\kappa}$?
\end{question}
\section{Successors of singulars}
So far, our results revolved around regular cardinals which are not the successors of singulars of cofinality below $\kappa$. In this section, we will see that the situation changes for such cardinals. First, let us recall some basic Boolean algebra terminology. Given a complete Boolean algebra $\mathbb{B}$, an ultrafilter $U$ on $\mathbb{B}$ is a subset of $\mathbb{B}$ such that $0_{\mathbb{B}}\notin U$ and $1_{\mathbb{B}}\in U$, it is $\leq_{\mathbb{B}}$-upwards closed, for every finite $\{b_0,...,b_{n-1}\}\subseteq U$, $b_0\wedge...\wedge b_{n-1}\in U$, and for every $b\in\mathbb{B}$, either $b$ or its complement is in $U$.
\begin{definition}
    We say that an ultrafilter $U$ on a complete Boolean algebra $\mathbb{B}$ is $\gamma$-indecomposable if for every maximal antichain $\{a_i\mid i<\gamma\}\subseteq \mathbb{B}$ there is $\delta<\gamma$ such that $\bigvee_{i<\delta}a_i\in U$.
\end{definition}
\begin{lemma}\label{lma:finest_partition}
    Suppose $\lambda$ is a limit cardinal,
    $\mathbb B$ is a $\lambda^+$-distributive complete Boolean algebra,
    and \(U\) is an ultrafilter on \(\mathbb B\) 
    that is \(\gamma\)-indecomposable for all \(\gamma\geq \lambda\). 
    Then there is a maximal antichain \(A\subseteq \mathbb B\) such that for any refinement \(B\) of \(A\),
    there is some \(u\in U\)
    such that \(B\mid u = A \mid u\).
    \begin{proof}
        Assume not. We construct a sequence $\langle A_\alpha : \alpha \leq \lambda\rangle$ of maximal antichains of \(\mathbb B\) by choosing \(A_{\alpha+1}\) so that there is no \(u\in U\) with \(A_{\alpha+1}\mid u = A_\alpha \mid u\), and
        taking least common refinements at limit stages using the $\lambda^+$-distributivity.

        By indecomposability, take \(u\in U\) such that \(|A_\lambda \mid u| < \lambda\). Let \(\gamma = |A_\lambda\mid u|\), and note that
        \(|A_\alpha\mid u|\leq \gamma\)
        for all \(\alpha\).

        If \(a\in A_{\gamma^+}\mid u\) and \(\alpha < \gamma^+\), denote by \(a_\alpha\) the unique \(\bar a \in A_\alpha\mid u\) such that \(a\leq \bar a\). 
        For distinct \(a,b\in A_{\gamma^+} \mid u\),
        let \(\alpha_{a,b} < \gamma^+\) be the least \(\alpha\) such that 
        \(a_\alpha \neq b_\alpha\). Note that
        \(\alpha_{a,b}\) exists since \(A_{\gamma^+}\mid u\) is the least common refinement of \(\langle A_\alpha\mid u  : {\alpha < \gamma^+}\rangle\).
        Let \(\beta < \gamma^+\) bound the \(\gamma\)-sized set of ordinals \(\alpha_{a,b}\). 
        Then \( A_\beta\mid u = A_{\gamma^+}\mid u\), which is a contradiction. 
    \end{proof}
\end{lemma}
\begin{lemma}\label{lma:no_indecomp_uf}
    Suppose $\lambda$ is a singular 
    limit of \(\lambda\)-strongly compact cardinals. 
    If $\mathbb B$ is a $\lambda^+$-distributive $\lambda^{++}$-cc atomless complete Boolean algebra that forces a new subset of $\lambda^+$,
    there are no \(\textnormal{cf}(\lambda)^+\)-complete \(\lambda^+\)-indecomposable ultrafilters on \(\mathbb B\).
    \begin{proof}
        Assume \(U\) is a counterexample.
        The Boolean ultrapower construction yields an elementary embedding \(j_U : V \to M_U\) and an \(M_U\)-generic filter \(G_U\subseteq j_U(\mathbb B)\).
        
        Applying Lemma \ref{lma:finest_partition}, fix an antichain \(A\subseteq \mathbb B\) such that for any refinement \(B\) of \(A\), there is some \(u\in U\) such that \(B \mid u = A \mid u\). By restricting to a measure one set, we can assume \(|A| < \lambda\). Let \(\eta < \lambda\) be a \(\lambda\)-strongly compact cardinal with \(\max\{|A|,\operatorname{cf}(\lambda)\}<\eta\).
        Let \(W\) be the projection of \(U\) onto \(P(A)\),
        so \(W\) is an ultrafilter over \(A\). 
        
        Let \(k : M_W\to M_U\) be the factor embedding, defined by \(k([f]_W) = j_U(f)(a)\) where \(a\) is the unique element of \(j_U(A)\cap G_U\). We claim that \(k\) is surjective.
        Each element of \(M_U\) has the form
        \(j_U(f)(b)\) where \(B\subseteq \mathbb B\) is a maximal antichain refining \(A\), \(f : B\to V\), and \(b\) is the unique element of \(j_U(B)\cap G_U\). Fix \(u\in U\) such that \(B \mid u = A \mid u\) and let \(\tilde f : A \to V\) be defined by \(\tilde f(x) = f(x')\) for the unique \(x'\in B\) such that \(x'\wedge u = x \wedge u\). Then \(k([\tilde{f}]_W) = j_U(\tilde{f})(a) = j_U(f)(b)\), as desired.

        Thus \(M_U = M_W\).
        By a slight generalization of Sakai's theorem \cite[Corollary 3.3]{Sakai} using the local approximation criterion \cite[Theorem~7.3.23]{GoldbergUA}, \(M_W\) has the \(\eta\)-approximation
        property in \(V\) for subsets of \(\lambda^+\). Take \(x\in M_W[G_U]\setminus M_W\) with \(x\subseteq j_W(\lambda^+)=\lambda^+\). By distributivity, \(x\) is \(\eta\)-approximated by \(M_W\).
        The local approximation property gives \(x\in M_W\), a contradiction. 
    \end{proof}
\end{lemma}

\begin{theorem}
    Suppose $\lambda$ is a singular 
    limit of \(\lambda\)-strongly compact cardinals
    and \(\mathcal U\) is a normal
    fine \(\kappa\)-complete ultrafilter on \(P_{\kappa}(\lambda^+)\) where
    \(\kappa > \textnormal{cf}(\lambda)\). Then Prikry forcing with \(\mathcal U\) does not project onto any $\lambda^+$-distributive atomless forcing. In particular, it does not project onto $\Add(\lambda^+,1)$.
    \begin{proof}
        Assume towards a contradiction that \(\mathbb B\) is an atomless \(\lambda^+\)-distributive projection of the supercompact Prikry forcing.
        Since the supercompact Prikry forcing is \(\lambda^+\)-centered, so is \(\mathbb B\).
        Thus \(\mathbb B\times\mathbb B\) is \(\lambda^{++}\)-cc, and by \cite[Note~1.3(2) and Lemma~1.5]{Usuba}, \(\mathbb B\) has the \(\lambda^{++}\)-approximation property.
        Since \(\mathbb B\) is atomless, forcing with \(\mathbb B\) adds a new set, and so by the \(\lambda^{++}\)-approximation property, it adds a new subset of \(\lambda^+\).

        Let \(j_{0,\omega} : V \to M_\omega\) be the \(\omega\)-th iterate of \(\mathcal U\). Then by Theorem \ref{Bukovsky-Dehornoy}, there is an \(M_\omega\)-generic ultrafilter \(G\subseteq j_{0,\omega}(\mathbb B)\)   such that $M_{\omega}[G]\subseteq \bigcap_{n<\omega}M_n$. By Lemma \ref{Central Lemma}, there is $n<\omega$ such that \[\sup\bigl(j_{0,\omega}(\lambda^+)\cap \text{Hull}^{M_\omega[G]}(j_{n,\omega}[M_n]\cup \{G\})\bigr)=\sup j_{n,\omega}''j_{n}(\lambda^+).\]
        Let \(U = j_{n,\omega}^{-1}[G]\). Since \(G\in\bigcap_{m<\omega}M_m\), we have \(U\in M_n\). Then \(U\)
        is a \(j_n(\kappa)\)-complete ultrafilter
        on \(\mathbb B^* = j_{0,n}(\mathbb B)\). Applying Lemma \ref{lma:no_indecomp_uf} in \(M_n\), we will obtain a contradiction by verifying that \(U\) is  \(j_{0,n}(\lambda^+)\)-indecomposable in \(M_n\).
        For this, suppose \(\langle a_\xi : \xi < j_{0,n}(\lambda^+)\rangle \in M_n\) enumerates a maximal antichain of \(\mathbb B^*\).
        We will show there is some \(\beta < j_{0,n}(\lambda^+)\) such that \(\bigvee_{\xi < \beta} a_\xi\in U\).
        Let \(\langle a^*_\xi : \xi < j_{0,\omega}(\lambda^+)\rangle = j_{n,\omega}(\langle a_\xi : \xi < j_{0,n}(\lambda^+)\rangle)\). Then the unique \(\zeta < j_{0,\omega}(\lambda^+)\) such that \(a^*_\zeta\in G\) belongs to \(\text{Hull}^{M_\omega[G]}(j_{n,\omega}[M_n]\cup \{G\})\).
        It follows that \(\zeta < \sup j_{n,\omega}''j_{0,n}(\lambda^+)\). So fix \(\beta < j_{0,n}(\lambda^+)\)
        such that \(\zeta < j_{n,\omega}(\beta)\). We then have
        \(\bigvee_{\xi < \beta} a_\xi \in U\)
        since \(j_{n,\omega}(\bigvee_{\xi < \beta} a_\xi)\in G\).
    \end{proof}
\end{theorem}

\textbf{Acknowledgment.} No generative AI was used in any way to produce the results of this paper. ChatGPT found the reference to Krueger's paper \cite{Krueger}. GPT-6.2 Sol was used to identify and correct typographical errors. The first author would like to thank Moti Gitik and Menachem Magidor for valuable comments and suggestions during the presentation of this work at the TAU Set Theory Seminar.
\bibliographystyle{alpha}
\bibliography{ref}

\begin{thebibliography}{BTW26}

\bibitem[AS83]{AbrahamShelah1994}
Uri Abraham and Saharon Shelah.
\newblock {F}orcing {C}losed {U}nbounded {S}ets.
\newblock {\em {T}he {J}ournal of {S}ymbolic {L}ogic}, 48(4):643--657, 1983.

\bibitem[BG21a]{partOne}
Tom Benhamou and Moti Gitik.
\newblock {I}ntermediate {M}odels of {M}agidor-{R}adin {F}orcing-{P}art {I}.
\newblock {\em Isr. J. of Mat.}, 2021.

\bibitem[BG21b]{TomMoti}
Tom Benhamou and Moti Gitik.
\newblock {S}ets in {P}rikry and {M}agidor {G}eneric {E}xtesions.
\newblock {\em Annals of Pure and Applied Logic}, 172(4):102926, 2021.

\bibitem[BG22]{Parttwo}
Tom Benhamou and Moti Gitik.
\newblock {I}ntermediate {M}odels of {M}agidor-{R}adin {F}orcing- part {I}{I}.
\newblock {\em Ann. Pure App. Log.}, 173:103107, 2022.

\bibitem[BG24]{BENHAMOU_GITIK_2024}
Tom Benhamou and Moti Gitik.
\newblock On {C}ohen and {P}rikry forcing notions.
\newblock {\em The Journal of Symbolic Logic}, 89(2):858–904, 2024.

\bibitem[BGH24]{Benhamou_Hayut_Gitik}
Tom Benhamou, Moti Gitik, and Yair Hayut.
\newblock The variety of projections of a tree prikry forcing.
\newblock {\em Journal of Mathematical Logic}, 24(03):2450006, 2024.

\bibitem[BTW26]{BenhamouTheiWeltsch2026}
Tom Benhamou, Sebastiano Thei, and Ben-Zion Weltsch.
\newblock On the intermediate models of strongly compact prikry forcing, 2026.

\bibitem[Dob03]{Dobrinen2003}
Natasha Dobrinen.
\newblock Games and general distributive laws in boolean algebras.
\newblock {\em Proceedings of the American Mathematical Society}, 131(1):309--318, 2003.

\bibitem[Git10]{Gitik2010}
Moti Gitik.
\newblock {\em {P}rikry {T}ype {F}orcings}, pages 1351--1447.
\newblock Springer Netherlands, Dordrecht, (2010).

\bibitem[GKK10]{PrikryCaseGitikKanKoe}
Moti Gitik, Vladimir Kanovei, and Peter Koepke.
\newblock {I}ntermediate {M}odels of {P}rikry {G}eneric {E}xtensions.
\newblock {\em http://www.math.tau.ac.il/~gitik/spr-kn.pdf}, page unpublished manuscript, 2010.

\bibitem[Gol22]{GoldbergUA}
G.~Goldberg.
\newblock {\em The Ultrapower Axiom}.
\newblock Berlin, Boston:De Gruyter, 2022.

\bibitem[Jec03]{Jech2003}
Thomas Jech.
\newblock {\em Set {T}heory}.
\newblock Springer Monographs in Mathematics. Springer-Verlag, Berlin, 2003.
\newblock The third millennium edition, revised and expanded.

\bibitem[KRS13]{MinimalPrikry}
Peter Koepke, Karen Rasch, and Philipp Schlicht.
\newblock {M}inimal {P}rikry-{T}ype {F}orcing for {S}ingularizing a {M}easurable {C}ardinal.
\newblock {\em J. Symb. Logic}, 78:85–--100, 2013.

\bibitem[Kru03]{Krueger}
John Krueger.
\newblock Fat sets and saturated ideals.
\newblock {\em The Journal of Symbolic Logic}, 68(3):837--845, 2003.

\bibitem[Lav78]{Laver1978}
Richard Laver.
\newblock Making the supercompactness of {$\kappa$} indestructible under {$\kappa$}-directed closed forcing.
\newblock {\em Israel Journal of Mathematics}, 29(4):385--388, 1978.

\bibitem[Mag77]{Magidor1977}
Menachem Magidor.
\newblock On the singular cardinals problem. i.
\newblock {\em Israel Journal of Mathematics}, 28(1-2):1--31, 1977.

\bibitem[Mag82]{Magidor1982}
Menachem Magidor.
\newblock Reflecting stationary sets.
\newblock {\em Journal of Symbolic Logic}, 47(4):755–771, 1982.

\bibitem[Mah47]{maharam}
Dorothy Maharam.
\newblock An {A}lgebraic {C}haracterization of {M}easure {A}lgebras.
\newblock {\em Annals of Mathematics}, 48:154–167, 1947.

\bibitem[Sac71]{Sacks1971}
Gerald~E. Sacks.
\newblock Forcing with perfect closed sets.
\newblock In {\em Axiomatic Set Theory}, volume~13 of {\em Proceedings of Symposia in Pure Mathematics}, pages 331--355. American Mathematical Society, Providence, RI, 1971.

\bibitem[Sak15]{Sakai}
Hiroshi Sakai.
\newblock Note on covering and approximation properties.
\newblock {\em RIMS K\^oky\^uroku}, 1949:81--88, 2015.

\bibitem[She03]{Shelah2003NotCollapsing}
Saharon Shelah.
\newblock Not collapsing cardinals $\leq\kappa$ in $(<\kappa)$-support iterations.
\newblock {\em Israel Journal of Mathematics}, 136:29--115, 2003.

\bibitem[Sol74]{Solovay1974}
Robert~M. Solovay.
\newblock Strongly compact cardinals and the {GCH}.
\newblock In {\em Proceedings of the Tarski Symposium}, volume~25 of {\em Proceedings of Symposia in Pure Mathematics}, pages 365--372. American Mathematical Society, Providence, RI, 1974.

\bibitem[Usu26]{Usuba}
Toshimichi Usuba.
\newblock Monotonicity of the ultrafilter number function.
\newblock {\em Archive for Mathematical Logic}, 65:545--568, 2026.

\end{thebibliography}
\end{document}